\documentclass[preprint, 3p]{elsarticle}

\usepackage{amssymb}
\usepackage{amsmath}
\usepackage{algorithm}
\usepackage{algpseudocode}
\usepackage{booktabs}
\usepackage{color}
\usepackage[font=small]{caption}
\usepackage{makecell}

\journal{Computer Methods in Applied Mechanics and Engineering}

\begin{document}

\begin{frontmatter}

\title{Wasserstein crossover for evolutionary algorithm-based topology optimization}

\author[osaka]{Taisei Kii}
\author[osaka]{Kentaro Yaji\corref{cor1}}
\ead{yaji@mech.eng.osaka-u.ac.jp}
\author[kansai]{Hiroshi Teramoto}
\author[osaka]{Kikuo Fujita}

\affiliation[osaka]{organization={Department of Mechanical Engineering, Graduate School of Engineering, The University of Osaka},
            addressline={2-1 Yamadaoka}, 
            city={Suita, Osaka},
            postcode={565-0871}, 
            country={Japan}}
\affiliation[kansai]{organization={Department of Mathematics, Faculty of Engineering Science, Kansai University},
            addressline={3-3-35 Yamate-cho}, 
            city={Suita, Osaka},
            postcode={564-8680}, 
            country={Japan}}
\cortext[cor1]{Corresponding author.}

\begin{abstract}
Evolutionary algorithms (EAs) are promising approaches for non-differentiable or strongly multimodal topology optimization problems, but they often suffer from the curse of dimensionality, generally leading to low-resolution optimized results.
This limitation stems in part from the difficulty of producing effective offspring through traditional crossover operators, which struggle to recombine complex parent design features in a meaningful way.
In this study, we propose a novel crossover operator for topology optimization, termed \textit{Wasserstein crossover}, and develop a corresponding EA-based optimization framework.
Our method leverages a morphing technique based on the Wasserstein distance---a distance metric between probability distributions derived from the optimal transport theory.
Its key idea is to treat material distributions as probability distributions and generate offspring as Wasserstein barycenters, enabling smooth and interpretable interpolation between parent designs while preserving their structural features.
The proposed framework incorporates Wasserstein crossover into an EA under a multifidelity design scheme, where low-fidelity optimized initial designs evolve through iterations of Wasserstein crossover and selection based on high-fidelity evaluation.
We apply the proposed framework to three topology optimization problems: maximum stress minimization in two- and three-dimensional structural mechanics, and turbulent heat transfer in two-dimensional thermofluids.
The results demonstrate that candidate solutions evolve iteratively toward high-performance designs through Wasserstein crossover, highlighting its potential as an effective crossover operator and validating the usefulness of the proposed framework for solving intractable topology optimization problems.
\end{abstract}




\end{frontmatter}

\section{Introduction}\label{sec1}

Topology optimization, first proposed by Bendsøe and Kikuchi~\cite{bendsoe1988}, offers a key advantage in its exceptionally high degree of design freedom.
Its fundamental concept is to optimize the material distribution within a given design domain by mathematical programming under a computational model of physical phenomena.
Typical topology optimization methods, such as the density-based method~\cite{bendsoe1989} and the level set method~\cite{sethian2000}, update design variables sequentially based on the sensitivity to the evaluation functions.
While these sensitivity-based approaches have high efficiency for convergence to reasonable optimized solutions, they inherently have significant challenges.
One of the major limitations is the need for differentiability of the evaluation functions, which is often not satisfied in practical applications.
For example, in stress-based topology optimization, the maximum stress is often approximated by $p$-norm or Kresselmeier-Steinhauser functions~\cite{le2010}, which are differentiable and continuous ones for the sake of sensitivity analysis.
Another challenge is the optimization of strongly multimodal functions.
A prominent example is turbulent flow in fluid topology optimization, where the complex physical phenomena lead to highly multimodal evaluation functions.
In such cases, oppressive parameter studies on optimization algorithms are often required to obtain a reasonable design, yet they still may be trapped in poor local optima.
As a result, most studies in fluid topology optimization have focused on laminar flows under low Reynolds number conditions~\cite{alexandersen2020}.
While such approximations have facilitated the application of topology optimization in a wide range of fields, there is potential to achieve even higher-performance designs through optimization approaches that can directly deal with the original evaluation functions and physics.

One of the promising approaches for solving such intractable optimization problems is the use of heuristic methods based on 
evolutionary algorithms (EAs)~\cite{thomas1996, yu2010}, which constitute a representative class of the broader field of evolutionary computation (EC).
Genetic algorithms (GAs)~\cite{goldberg1989}, which are population-based methods inspired by the process of natural selection, are a typical example of EAs.
In GAs, candidate solutions are encoded as chromosomes and handled as individuals within a population.
The population evolves over successive generations through the repeated three genetic operators: selection, which favors individuals with higher fitness; crossover, which recombines features from selected parents to generate offspring; and mutation, which introduces random variations.
In these optimization processes, the objective and constraint functions are evaluated as fitness values, without requiring their sensitivity information.
Furthermore, compared to the local search behavior of sensitivity-based methods, the population-based emergent optimization mechanisms of EAs enable effective global search.
As a result, they are well-suited for tackling intractable optimization problems that are non-differentiable or exhibit strong multimodality.

Recently, various EA-based topology optimization methods have been proposed, with innovations in both encoding material distributions as chromosomes and selecting appropriate algorithms, as summarized by Guirguis et al.~\cite{guirguis2020}.
The simplest approach, as developed by Chapman et al.~\cite{chapman1994}, Jakiela et al.~\cite{jakiela2000}, and Wang and Tai~\cite{wang2005}, is GA-based topology optimization with bit array representation.
As examples of GA-based approaches employing different encoding schemes, Madeira et al.~\cite{madeira2010} and Nimura and Oyama~\cite{nimura2024} proposed graph and quadtree representation-based methods, respectively.
As other examples using different algorithms within the broader EC framework with bit array representation, Wu and Theng~\cite{wu2010} and Luh et al.~\cite{luh2011} employed differential evolution and particle swarm optimization, respectively.
While these methods can also be applied to intractable problems due to the strength of EAs, EA-based topology optimization is typically constrained by a limitation commonly known as the curse of dimensionality, which refers to the exponential deterioration of search performance as the number of design variables increases.
Due to this limitation, most EAs have been applied to optimization problems with relatively few design variables, typically on the order of a thousand or less, and EA-based topology optimization has likewise been carried out on coarsely discretized design domains.
Sigmund~\cite{sigmund2011} has emphasized that such coarse discretizations, imposed by the curse of dimensionality, inevitably lead to low-resolution optimized designs, and highlighted this limitation contrasting with the high-resolution designs achievable by the density-based method.
Nevertheless, several pioneering works have sought to overcome this challenge by employing more sophisticated encoding schemes or advanced algorithms to achieve finer optimized designs.
For example, Fujii~\cite{fujii2018} and Tanaka and Fujii~\cite{tanaka2024} developed approaches that combine level set boundary representations with an advanced EA known as the covariance matrix adaptation evolution strategy (CMA-ES).
In addition, Furuta et al.~\cite{furuta2024} incorporated a Karhunen-Lo{\`e}ve expansion, which enables more compact representations of material distributions, into a differential evolution algorithm.

As discussed above, although a wide variety of EA-based topology optimization methods have been developed, most of them focus on how to represent material distributions as chromosomes and which algorithm to use.
While these algorithmic aspects are important in EA-based approaches, it is also crucial to consider how effectively new candidate solutions can be generated within the algorithm.
Especially, as crossover plays its central role in EA-based approaches, it is one of the key factors that largely determines the overall optimization performance of EAs.
For example, the simplest crossover strategies, such as exchanging parts of material distributions represented as $\{0, 1\}$ binary vectors, or linearly interpolating ones represented as $[0, 1]$ continuous vectors in the manner of the density-based method, often lead to disconnected designs.
In other words, in topology optimization, traditional crossover operators struggle to recombine parent design features in a meaningful way.
Therefore, we believe that the design of effective crossover operators is essential for the successful development of EA-based topology optimization.

As one of the pioneering works aiming to overcome the limitations in EA-based topology optimization through the development of crossover strategies, Yamasaki et al. proposed a framework called data-driven topology design (DDTD)~\cite{yamasaki2021}.
Its basic concept is to combine EAs with the idea of data-driven design~\cite{woldseth2022} by employing a deep generative model to perform the role of crossover.
Deep generative models, such as variational autoencoders (VAEs)~\cite{kingma2013} and generative adversarial networks (GANs)~\cite{goodfellow2014}, are capable of generating an arbitrary dataset similar to those in the training dataset.
DDTD takes advantage of this property by employing the VAE as a crossover operator to produce new candidate solutions that inherit features from a parent dataset, and iteratively updates the population based on the GA.
Leveraging the strength of EAs for intractable problems, DDTD has recently been applied to wide range of topology optimization problems such as structural mechanics~\cite{yamasaki2021, kii2024, kato2025, yang2025pca, yang2025ifl}, thermofluids~\cite{yang2025ifl, yaji2022, luo2025a, luo2025b, urata2025}, and electromagnetics~\cite{zhou2025}, including more advanced applications like solid-porus infill structure design~\cite{xu2025} and concurrent optimization of multiple design variable fields~\cite{kawabe2025}.
While the use of VAEs serves as a rational crossover operator for topology optimization, their inherent nature poses certain challenges.
Specifically, the latent space constructed by compressing the material distributions of the population typically corresponds to a limited subspace of the vast overall design space in topology optimization.
Consequently, generating new candidates through sampling from this subspace raises concerns regarding its suitability as an emergent crossover mechanism capable of global search.
Moreover, it is also known that VAEs have limitations in the data dimensionality which can be effectively learned~\cite{yang2025pca, yang2025ifl}, which requires incorporating auxiliary techniques such as principal component analysis (PCA)~\cite{yang2025pca} or image fragmented learning~\cite{yang2025ifl} so that DDTD can be applied to three-dimensional topology optimization problems.
Thus, a crossover operator based on generative models does not fully resolve the challenges in EA-based topology optimization, and the development of more effective crossover operators remains an important issue.

In this study, we introduce a novel crossover operator for topology optimization, termed \textit{Wasserstein crossover}, which applies a morphing technique based on the Wasserstein distance~\cite{solomon2015}.
The Wasserstein distance~\cite{villani2009}, which is also known as the earth mover's distance~\cite{rubner2000}, is a distance metric between probability distributions based on optimal transport theory~\cite{peyre2019}.
Solomon et al.~\cite{solomon2015} developed a morphing technique for its application to graphics.
This morphing approach involves treating pixel or voxel images as probability distributions and calculating their intermediate distribution based on optimal transport, i.e., the Wasserstein barycenter.
Inspired by this idea, we adopt density-based representations of material distributions, following the density-based topology optimization method, and generate offspring by computing Wasserstein barycenters of the parent distributions.
Conceptually, this amounts to interpolating between parent material distributions, yielding a simple yet distinctive crossover operator for EA-based topology optimization.
Owing to the unique interpolative behavior induced by optimal transport, Wasserstein crossover can effectively recombine parent design features and generate diverse candidate solutions, which are particularly advantageous for EAs.
We propose a novel EA-based topology optimization framework incorporating Wasserstein crossover into DDTD under the multifidelity design scheme~\cite{yaji2020}, where initial designs derived by solving a low-fidelity problem are updated through iterative cycles of high-fidelity evaluation, selection, and Wasserstein crossover.
This paper applies the proposed framework to several topology optimization problems, which are known to be intractable in structural mechanics and thermofluids, including a three-dimensional case, and demonstrates its effectiveness.

The rest of this paper is organized as follows.
Section~\ref{sec2} provides background on the Wasserstein distance derived from optimal transport theory and its application to morphing, highlighting its potential as a crossover operator through comparison with other morphing approaches.
Section~\ref{sec3} describes Wasserstein crossover and the proposed framework incorporating it into DDTD.
Section~\ref{sec4} presents the detailed numerical implementation of the proposed framework, including the overall algorithm and Wasserstein crossover.
Section~\ref{sec5} demonstrates the application of the proposed framework to several intractable topology optimization problems and discusses the results.
Finally, Section~\ref{sec6} concludes the paper.

\section{Preliminaries}\label{sec2}

In this section, we provide some prior knowledge on the Wasserstein distance and its application to graphical morphing.
Subsequently, we compare it with other morphing methods, including a VAE-based one employed in DDTD, to highlight its potential as a crossover operator for topology optimization.
We also review related works that have applied the idea of the Wasserstein distance to topology optimization and clarify the focus of this paper by highlighting the differences from those prior works.

\subsection{Wasserstein distance}\label{sec21}

The fundamental idea of optimal transport is to measure the minimal cost required to transform one probability distribution $\mu$ into another distribution $\nu$ defined on $\mathbb{R}^d$ (typically with $d=1, 2, 3$)~\cite{villani2009}, which is mathematically formulated as follows:
\begin{equation}\label{eq_ot}
  C(\mu, \nu)=\inf_{\pi \in \varPi(\mu, \nu)} \int_{\mathbb{R}^d \times \mathbb{R}^d} c(\mathbf{x}, \mathbf{y})\,\mathrm{d}\pi(\mathbf{x}, \mathbf{y}),
\end{equation}
where $c(\mathbf{x}, \mathbf{y})$ denotes the cost of transporting a unit mass from $\mathbf{x}\in\mathbb{R}^d$ to $\mathbf{y}\in\mathbb{R}^d$, and $\varPi(\mu, \nu)$ is the set of all transport plans with marginals $\mu$ and $\nu$.
The Wasserstein distance $W_p(\mu, \nu)$ is defined as $c(\mathbf{x}, \mathbf{y})$ with a cost function based on the $L^p$ norm $\|\mathbf{x} - \mathbf{y}\|_p$, where $p\geq 1$, to satisfy the  the axioms of a distance function, as follows:
\begin{equation}\label{eq_wasserstein}
  W_p(\mu, \nu) = \left( \inf_{\pi \in \varPi(\mu, \nu)} \int_{\mathbb{R}^d \times \mathbb{R}^d} \|\mathbf{x} - \mathbf{y}\|_p^p\,\mathrm{d}\pi(\mathbf{x}, \mathbf{y})\right)^{1/p}.
\end{equation}
$W_p$ in Eq.~\eqref{eq_wasserstein} is the Wasserstein distance of order $p$, or simply the $p$-Wasserstein distance.
In particular, the case of $p=1$ is also known as the earth mover's distance~\cite{rubner2000}.
By explicitly incorporating the ground cost $\|\mathbf{x} - \mathbf{y}\|_p^p$ into its definition, the Wasserstein distance reflects the geometric structure of the underlying space, providing an intuitive interpretation as the minimal effort required to morph one distribution into another.
Unlike the simple measures such as the Kullback-Leibler divergence, which can become undefined or infinite when the distributions do not overlap, the Wasserstein distance remains well-defined and captures differences both in location and shape of the distributions.
These properties make it a particularly expressive and versatile metric for comparing probability distributions.

While the Wasserstein distance is such a powerful metric, it is also known to be computationally expensive.
To address this issue, a regularized formulation is often used by augmenting the objective function in Eq.~\eqref{eq_wasserstein} with an additional regularization term~\cite{cuturi2013}, as follows:
\begin{equation}\label{eq_wasserstein_reg}
  W_{p, \varepsilon}(\mu, \nu) = \left( \inf_{\pi \in \varPi(\mu, \nu)} \int_{\mathbb{R}^d \times \mathbb{R}^d} \|\mathbf{x} - \mathbf{y}\|_p^p\,\mathrm{d}\pi(\mathbf{x}, \mathbf{y}) - \varepsilon H(\pi)\right)^{1/p},
\end{equation}
where $\varepsilon > 0$ is a regularization coefficient, and $H(\pi)$ is the entropy of the transport plan $\pi$, defined as follows:
\begin{equation}\label{eq_entropy}
  H(\pi) = -\int_{\mathbb{R}^d \times \mathbb{R}^d} \left( \log\pi(\mathbf{x}, \mathbf{y})-1 \right)\,\mathrm{d}\pi(\mathbf{x}, \mathbf{y}).
\end{equation}
$W_{p, \varepsilon}$ in Eq.~\eqref{eq_wasserstein_reg} is known as the entropy-regularized Wasserstein distance, and one of its advantages is that the regularization term $H(\pi)$ renders the objective function strongly convex, which ensures smooth convergence to the optimal solution.

Another key advantage of the entropic regularization is that it enables much faster computation of the approximated optimal solution using the Sinkhorn algorithm~\cite{cuturi2013, sinkhorn1967}.
To describe this, we consider the discretized case of the entropy-regularized Wasserstein distance, where the input probability distributions $\mu$ and $\nu$ are represented as discrete vectors $\mathbf{a}, \mathbf{b}\in\mathbb{R}^n$, and the distance between discrete points is given by a cost matrix $\mathbf{C}\in\mathbb{R}^{n\times n}$.
The regularized optimal transport problem corresponding to the continuous case in Eq.~\eqref{eq_wasserstein_reg} can then be formulated as follows:
\begin{eqnarray}\label{eq_sinkhorn}
  \begin{aligned}
    & \underset{\mathbf{P}\in\mathbb{R}^{n\times n}}{\text{minimize}} && \sum_{i,j}P_{ij}C_{ij}+\varepsilon\sum_{i,j}P_{ij} (\log P_{ij}-1) \\
    & \text{subject to} && \mathbf{P}\mathbf{1} = \mathbf{a}, \\
    &&& \mathbf{P}^\top\mathbf{1} = \mathbf{b}, \\
    &&& P_{ij} \geq 0,
  \end{aligned}
\end{eqnarray}
where $\mathbf{P}$ is the transport plan corresponding to $\pi$ in the continuous formulation of Eq.~\eqref{eq_wasserstein_reg}, and $\mathbf{1}\in\mathbb{R}^n$ is a vector of ones.
For the optimization problem in Eq.~\eqref{eq_sinkhorn}, the pseudocode of the Sinkhorn algorithm, derived from its Lagrangian dual formulation, is presented in Algorithm~\ref{alg_sinkhorn}.
In this way, the Sinkhorn algorithm solves the entropy-regularized optimal transport problem by iteratively updating the scaling vectors $\mathbf{u}$ and $\mathbf{v}$, with a computational complexity of $\mathcal{O}(n^2)$ per iteration.
In contrast, the original optimal transport problem formulated without the regularization term in Eq.~\eqref{eq_sinkhorn} becomes a linear programming problem, whose computational complexity is $\mathcal{O}(n^3)$.
Although the number of Sinkhorn iterations in Algorithm~\ref{alg_sinkhorn} depends on the strength of the regularization, in practice, convergence is achieved within tens or hundreds of iterations.
Thus, the combination of the entropic regularization and Sinkhorn algorithm serves as an efficient approximation of the optimal transport solution with significantly reduced computational cost.
Furthermore, as shown in Algorithm~\ref{alg_sinkhorn}, the Sinkhorn algorithm essentially consists of simple matrix-vector computations, which makes it well-suited for GPU acceleration, as demonstrated by Ryu et al.~\cite{Ryu2018}.
In particular, when the cost function $c(\mathbf{x}, \mathbf{y})$ is given by the squared Euclidean norm $\|\mathbf{x} - \mathbf{y}\|_2^2$, which corresponds to the case of the 2-Wasserstein distance $W_{2, \varepsilon}$, the matrix $\mathbf{K}$ in Algorithm~\ref{alg_sinkhorn} becomes a Gaussian kernel matrix.
As a result, the computations of $\mathbf{K} \mathbf{v}$ and $\mathbf{K}^\top \mathbf{u}$ are equivalent to convolutions with a Gaussian filter, which enables further acceleration through efficient convolution-based implementations.
Thus, the Wasserstein distance, defined as a metric between probability distributions based on optimal transport, provides a theoretically grounded measure for which efficient computational methods have also been well established.

\begin{algorithm}[ht]
  \caption{Sinkhorn algorithm for computing the entropy-regularized optimal transport}\label{alg_sinkhorn}
  \begin{algorithmic}[1]
    \Require Probability vectors $\mathbf{a}, \mathbf{b}\in\mathbb{R}^n$, cost matrix $\mathbf{C}\in\mathbb{R}^{n\times n}$, entropic regularization coefficient $\varepsilon > 0$
    \State $\mathbf{K} \gets \exp\left(-{\mathbf{C}}/{\varepsilon}\right)$ \Comment{Element-wise exponential}
    \State $\mathbf{u} \gets \mathbf{1}$, $\mathbf{v} \gets \mathbf{1}$ \Comment{Initialize with ones}
    \Repeat \Comment{Sinkhorn iterations}
      \State $\mathbf{u} \gets \mathbf{a} \oslash (\mathbf{K} \mathbf{v})$ \Comment{$\oslash$ denotes element-wise division}
      \State $\mathbf{v} \gets \mathbf{b} \oslash (\mathbf{K}^\top \mathbf{u})$
      \State $r_a \gets \|\mathbf{u} \odot (\mathbf{K}\mathbf{v}) - \mathbf{a}\|_1 \quad r_b \gets \|\mathbf{v} \odot (\mathbf{K}^\top \mathbf{u}) - \mathbf{b}\|_1$ \Comment{Marginal residuals, $\odot$ denotes element-wise multiplication}
    \Until{$\max(r_a, r_b) < \tau$} \Comment{Convergence criterion with tolerance $\tau$}
    \State $\mathbf{P} = \operatorname{diag}(\mathbf{u}) \mathbf{K} \operatorname{diag}(\mathbf{v})$ \Comment{Transport plan}
    \State \Return $\sum_{i,j} P_{ij} C_{ij}$ \Comment{Wasserstein distance}
  \end{algorithmic}
\end{algorithm}

\subsection{Wasserstein barycentric morphing}\label{sec22}

The concepts of the Wasserstein distance and optimal transport have been applied in various fields, including deep learning~\cite{arjovsky2017}, image processing~\cite{rubner2000}, and natural language processing~\cite{kusner2015}.
Among these, one notable application is morphing proposed by Solomon et al.~\cite{solomon2015}, where these techniques are used to interpolate between different shapes or distributions.
It is based on the concept of the Wasserstein barycenter~\cite{agueh2011}, which is defined as the weighted average of multiple probability distributions concerning the Wasserstein distance, given as follows:
\begin{equation}\label{eq_barycenter}
  \mu ^* = \arg\min_{\mu} \sum_{i=1}^N \lambda_i W_p(\mu, \mu_i),
\end{equation}
where $\mu_i$ are the $N$ input probability distributions, $\lambda_i$ are the weights satisfying $\sum_{i=1}^N \lambda_i = 1$.
The Wasserstein barycenter in Eq.~\eqref{eq_barycenter} is formulated as a minimization problem for the multiple Wasserstein distances, which themselves require solving optimization problems in Eq.~\eqref{eq_wasserstein}.
This suggests that computing the Wasserstein barycenter would require solving a nested or two-level optimization problem.
However, similar to the regularized Wasserstein distance in Eq.~\eqref{eq_wasserstein_reg}, entropic regularization can also be considered for the barycenter computation, allowing for an efficient approximation by a matrix-based iterative algorithm proposed by Benamou et al.~\cite{benamou2015}, as outlined in Algorithm~\ref{alg_barycenter}.
Owing to the regularization, the approximated Wasserstein barycenter can be obtained without solving complex nested optimization problems, instead relying on simple iterative matrix operations.
In addition, as desctibed in Section~\ref{sec21}, these matrix operations can be further accelerated using GPU implementations, and in the case of the 2-Wasserstein distance, convolution with a Gaussian filter can also be used to speed up the computaion of $\mathbf{K}\mathbf{v}^{(i)}$ and $\mathbf{K}^\top\mathbf{u}^{(i)}$ in Algorithm~\ref{alg_barycenter}.
Note that the entropic regularization coefficient $\varepsilon$ and the convergence tolerance $\tau$ in Algorithm~\ref{alg_barycenter} affect the number of Sinkhorn iterations and should be chosen to balance computational cost and the quality of the approximated barycenter.

\begin{algorithm}[ht]
  \caption{Sinkhorn algorithm for computing the Wasserstein barycenter}\label{alg_barycenter}
  \begin{algorithmic}[1]
    \Require Probability vectors $\mathbf{a}^{(1)}, \ldots, \mathbf{a}^{(N)}\in\mathbb{R}^n$, weights $\lambda_1, \ldots, \lambda_N$, cost matrix $\mathbf{C}\in\mathbb{R}^{n\times n}$, entropic regularization coefficient $\varepsilon > 0$
    \State $\mathbf{K} \gets \exp\left(-{\mathbf{C}}/{\varepsilon}\right)$ \Comment{Element-wise exponential}
    \State $\mathbf{u}^{(i)} \gets \mathbf{1}$, $\mathbf{v}^{(i)} \gets \mathbf{1} \quad \left(\forall i = 1, \ldots, N\right)$ \Comment{Initialize with ones}
    \Repeat \Comment{Sinkhorn iterations}
      \For{$i = 1, \ldots, N$}
        \State $\mathbf{u}^{(i)} \gets \mathbf{a}^{(i)} \oslash \left( \mathbf{K} \mathbf{v}^{(i)}\right)$
        \State $\mathbf{v}^{(i)} \gets \left(\prod_j \left(\mathbf{K}^\top \mathbf{u}^{(j)}\right)^{\lambda_j}\right) \oslash \left( \mathbf{K}^\top \mathbf{u}^{(i)}\right)$
      \EndFor
      \State $E \gets \sum_{e=1}^n \mathrm{Std}\big(\{\mathbf{v}^{(i)} \odot (\mathbf{K}^\top \mathbf{u}^{(i)})\}_{i=1}^N\big)_e$ \Comment{Marginal deviations}
    \Until{$E < \tau$} \Comment{Convergence tolerance $\tau$}
    \State \Return $\mathbf{a}^{*} = \prod_i \left( \mathbf{K}^\top \mathbf{u}^{(i)}\right)^{\lambda_i}$ \Comment{Compute barycenter}
  \end{algorithmic}
\end{algorithm}

Any shape data represented in a pixel or voxel grid can be treated as a probability distribution by normalizing its sum to one.
By gradually varying the weight $\lambda_i$ in Eq.~\eqref{eq_barycenter} for each input distribution $\mu_i$, the Wasserstein barycenters can be computed to interpolate between the original shapes, resulting in smooth morphing based on optimal transport.
Benamou et al.~\cite{benamou2015} demonstrated this morphing technique using three 2D shapes (diamond, ring, and square), and Solomon et al.~\cite{solomon2015} extended it to 3D shapes.

\subsection{Comparison with other morphing methods}\label{sec23}

\begin{figure}[t]
  \centering
  \includegraphics[width=\textwidth]{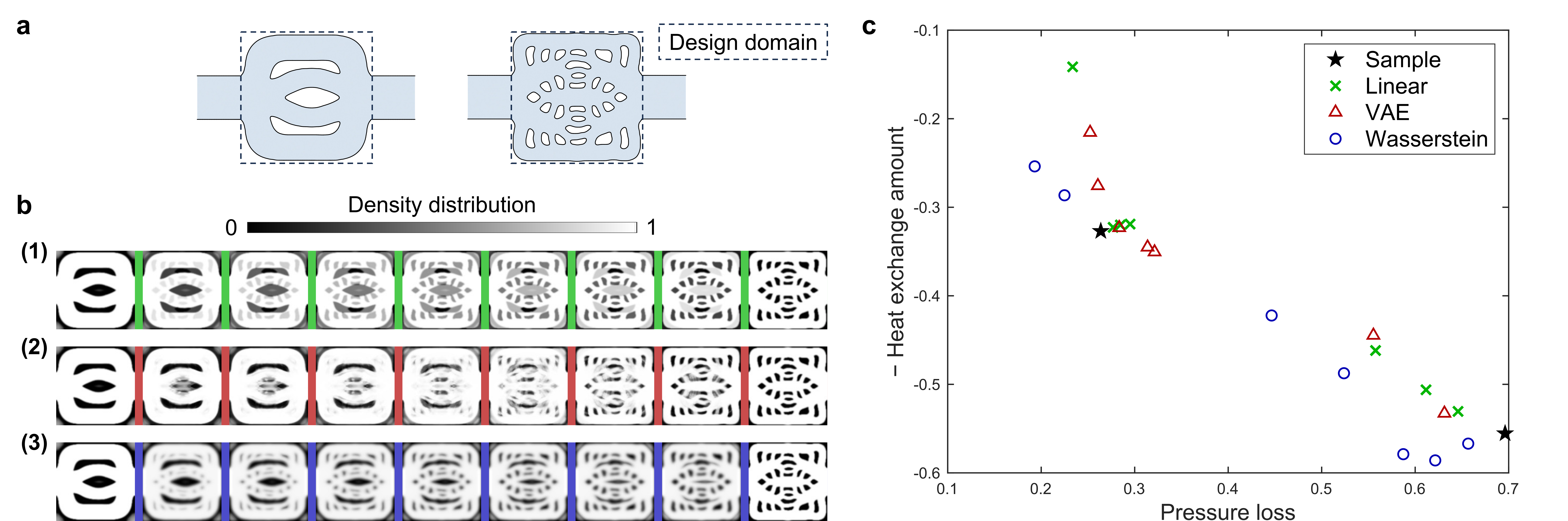} 
  \caption{Comparison of different morphing methods between two sample data for a heat transfer problem: (a) two sample flow channels; (b) morphing results for design domain using (1) linear interpolation, (2) a variational autoencoder (VAE), and (3) Wasserstein barycentric interpolation; (c) objective values (pressure loss versus negative heat exchange amount) of sample and interpolated flow channels.}\label{fig_morphing_compare}
\end{figure}

We compare the Wasserstein barycentric morphing described in Section~\ref{sec22} with two alternative morphing methods to examine their interpolation behaviors and verify their potential as a crossover operator for topology optimization.
Figure~\ref{fig_morphing_compare} shows the interpolated results obtained by morphing two different sample flow channels for the heat transfer problem in Section~\ref{sec52}, together with their performance values.
As comparison methods, we consider:
\begin{enumerate}\def\labelenumi{(\theenumi)}
  \item linear interpolation, the simplest approach, which corresponds to Euclidean barycentric interpolation;
  \item VAE-based interpolation, obtained by linearly interpolating between two latent vectors.
  \item Wasserstein barycentric interpolation, proposed in this paper.
\end{enumerate}
Comparing the three approaches in Fig.~\ref{fig_morphing_compare}(b), linear interpolation results in a simple superposition of the two channels, whereas VAE-based interpolation exhibits a more complex and gradual transition through their intermediate shapes. 
This behavior stems from the nature of VAEs: the encoder and decoder networks are trained on a certain dataset, and interpolations in the constructed latent space implicitly incorporate the influence of the entire training dataset---100 different flow channels in this case---even when two specific designs are interpolated.
As a result, the generated intermediate shapes reflect not only the two sample designs but also the characteristics of other flow channels in the training dataset.
In contrast, the Wasserstein barycentric interpolation is derived from explicitly considering optimal transport between two sample designs.
Compared to the grayscale-based transitions in the other two approaches, the Wasserstein barycentric morphing yields qualitatively different outcomes in that the fluid and solid regions (represented in white and black, respectively) remain relatively well distinguished throughout the transition.
The resulting morphing involves a physically meaningful interpolation, in which the topology and branching patterns of flow channels are gradually transitioned.

Focusing on the performance values of the sample and interpolated channels shown in Fig.~\ref{fig_morphing_compare}(c), we observe that the linear interpolated designs tend to exhibit mostly comparable performance values to those of the original samples.
In contrast, the VAE-based interpolations are more scattered around the sample points, indicating slightly greater variation.
Notably, however, the Wasserstein barycentric interpolations are more widespread in the objective space, including some designs that dominate the original samples in terms of Pareto optimality.
These high-performance designs can be attributed to the distinctive interpolation behavior afforded by optimal transport, where the morphing in Fig.~\ref{fig_morphing_compare}(b3) yields hybrid flow channels that effectively recombine the features of the samples.

These results demonstrate the potential of the Wasserstein barycentric morphing as a crossover operator, in that it enables the generation of offspring by shape interpolation between parents, some of which can have better performance compared to their parents.
While the examples presented in Fig.~\ref{fig_morphing_compare} focus on two-dimensional flow channels for a heat transfer problem, we further investigate the performance of the proposed Wasserstein crossover through several numerical examples, including structural mechanics problems and three-dimensional cases, in Section~\ref{sec5}.

\subsection{Novelty in this paper compared to related works}\label{sec24}

A representative prior work that incorporates the concept of optimal transport and Wasserstein distance into topology optimization is the use of deep generative models based on the so-called Wasserstein GAN~\cite{arjovsky2017}.
Wasserstein GANs are a variant of generative adversarial networks in which the earth mover’s distance, i.e., the 1-Wasserstein distance, is used as the loss function.
This formulation is known to improve learning stability and enhance the quality of the generated images.
Recently, several topology optimization methods using Wasserstein GANs have been developed~\cite{rawat2019a, rawat2019b, wang2023, pereira2024, zeng2024}, but most of them aim to reduce computational cost by learning from a pre-topology-optimized dataset to predict optimized designs and eliminate iterations in topology optimization.
This objective differs fundamentally from the aim of this paper, which focuses on developing a crossover operator within an EA-based topology optimization framework.

As an example of incorporating the idea of optimal transport into EAs, Sow et al.~\cite{sow2024} proposed crossover and mutation operators based on Wasserstein barycenters.
While their work focuses on the optimization of point clouds and considers wind-farm layout design as an engineering application, the design variables in their setting consist of the positions of only a few dozen points.
This is fundamentally different from topology optimization, where the design variables typically involve material distributions over grids with thousands or more degrees of freedom.

Another example of applying Wasserstein barycenters to shape design is the framework proposed by Ma et al.~\cite{ma2021} for the design of underwater swimming robots.
They focused on the co-design of both the geometry and controller of swimmers, and successfully utilized Wasserstein barycentric interpolation to generate candidate designs of robots that are represented with tens of thousands of voxel resolutions.
However, their approach remains centered on repeated simulation and interpolation between base designs given by the user, which is fundamentally different from the idea of EAs, where a population of candidate designs is evolved through genetic operations.

In this paper, we propose an EA-based topology optimization framework that introduces Wasserstein barycentric interpolation as an effective crossover operator for material distributions.
The originality and novelty of this paper, compared with those related works, can be summarized as follows:
\begin{itemize}
  \item Introduction of Wasserstein crossover as a dedicated and effective crossover operator for material distributions, based on the Wasserstein barycentric morphing.
  \item Development of an EA-based topology optimization framework capable of handling finely discretized continuum representations of material distributions with reasonable computational cost.
\end{itemize}

\section{Framework}\label{sec3}

In this section, we describe the proposed framework incorporating Wasserstein crossover.
Our proposed framework is based on DDTD~\cite{yamasaki2021, yaji2022}, and focuses on solving the general multi-objective topology optimization problem, which is formulated as follows:
\begin{equation}\label{eq_hf}
  \begin{aligned}
    & \underset{\gamma (\mathbf{x})}{\text{minimize}} && [J_1(\gamma), J_2(\gamma) \ldots , J_{n_{\text{o}}}(\gamma)], \\
    & \text{subject to} && G_j(\gamma) \leq 0 \quad (j=1, 2, \ldots , n_{\text{c}}),\\
    &&& \gamma (\mathbf{x})\in\{0, 1\}, \quad \forall\mathbf{x}\in D,\\
  \end{aligned}
\end{equation}
where $\gamma (\mathbf{x})$ is the binary material distribution defined on the design domain $D\subset\mathbb{R}^d$ ($d=2\text{ or }3$), $J_i(\gamma)$ and $G_j(\gamma)$ are the objective and constraint functions, respectively.
It is often difficult to directly solve the optimization problem of Eq.~\eqref{eq_hf} due to the non-differentiability and strong multimodality of evaluation functions $J_i(\gamma)$ and $G_j(\gamma)$.
To address this issue, we employ a multifidelity design approach~\cite{yaji2020}, treating the formulation of Eq.~\eqref{eq_hf} as a high-fidelity (HF) model and constructing a low-fidelity (LF) model, which is formulated as a simplified pseudo-problem as follows:
\begin{equation}\label{eq_lf}
  \begin{aligned}
    & \underset{\boldsymbol{\gamma}^{(k)}}{\text{minimize}} && \widetilde{J}(\boldsymbol{\gamma}^{(k)}, \mathbf{s}^{(k)}), \\
    & \text{subject to} && \widetilde{G}_j(\boldsymbol{\gamma}^{(k)}, \mathbf{s}^{(k)}) \leq 0 \quad (j=1, 2, \ldots , \widetilde{n}_{\text{c}}), \\
    &&& \boldsymbol{\gamma}^{(k)}\in [0, 1]^n, \\
    & \text{for given} && \mathbf{s}^{(k)},
  \end{aligned}
\end{equation}
where $\widetilde{J}(\boldsymbol{\gamma}^{(k)}, \mathbf{s}^{(k)})$ and $\widetilde{G}_j(\boldsymbol{\gamma}^{(k)}, \mathbf{s}^{(k)})$ are the objective and constraint functions of the LF model, respectively.
Solving the optimization problem of Eq.~\eqref{eq_lf} is referred to as LF optimization, whose primary aim is to derive candidate solutions for the original problem of Eq.~\eqref{eq_hf}.
The vector composed of artificial parameters $\mathbf{s}=[s_1, s_2, \ldots, s_{n_\text{sd}}]^\top$ is called a seeding parameter, and $\mathbf{s}^{(k)}$ with $k=1, 2, \ldots, N_{\text{lf}}$ represents the sample point of $\mathbf{s}$, which serves to induce diverse candidate designs.
Taking into account the compatibility with Wasserstein crossover, which interpolates between candidate solutions represented as probability distributions, the LF optimization problem is assumed to be solved using the density-based method~\cite{bendsoe1989}.
Thus, the objective and constraint functions, $\widetilde{J}(\boldsymbol{\gamma}^{(k)}, \mathbf{s}^{(k)})$ and $\widetilde{G}_j(\boldsymbol{\gamma}^{(k)}, \mathbf{s}^{(k)})$, are assumed to be differentiable with respect to the design variable $\boldsymbol{\gamma}^{(k)}$, and the LF optimization problem of Eq.~\eqref{eq_lf} should be easily solvable by gradient-based optimizers.
Accordingly, the design variable $\boldsymbol{\gamma}^{(k)}$ is defined as a continuous vector representing a material density distribution with values in the range $[0, 1]$, in contrast to the original binary material distribution $\gamma (\mathbf{x})$ defined in the continuous domain.
In addition, the LF optimization problem of Eq.~\eqref{eq_lf} is converted to a single-objective problem from the original multi-objective problem of Eq.~\eqref{eq_hf} using the weighted sum method~\cite{zadeh1963}, $\varepsilon$-constraint method~\cite{haimes1971}, and so on.

\begin{figure}[t]
  \centering
  \includegraphics[width=\textwidth]{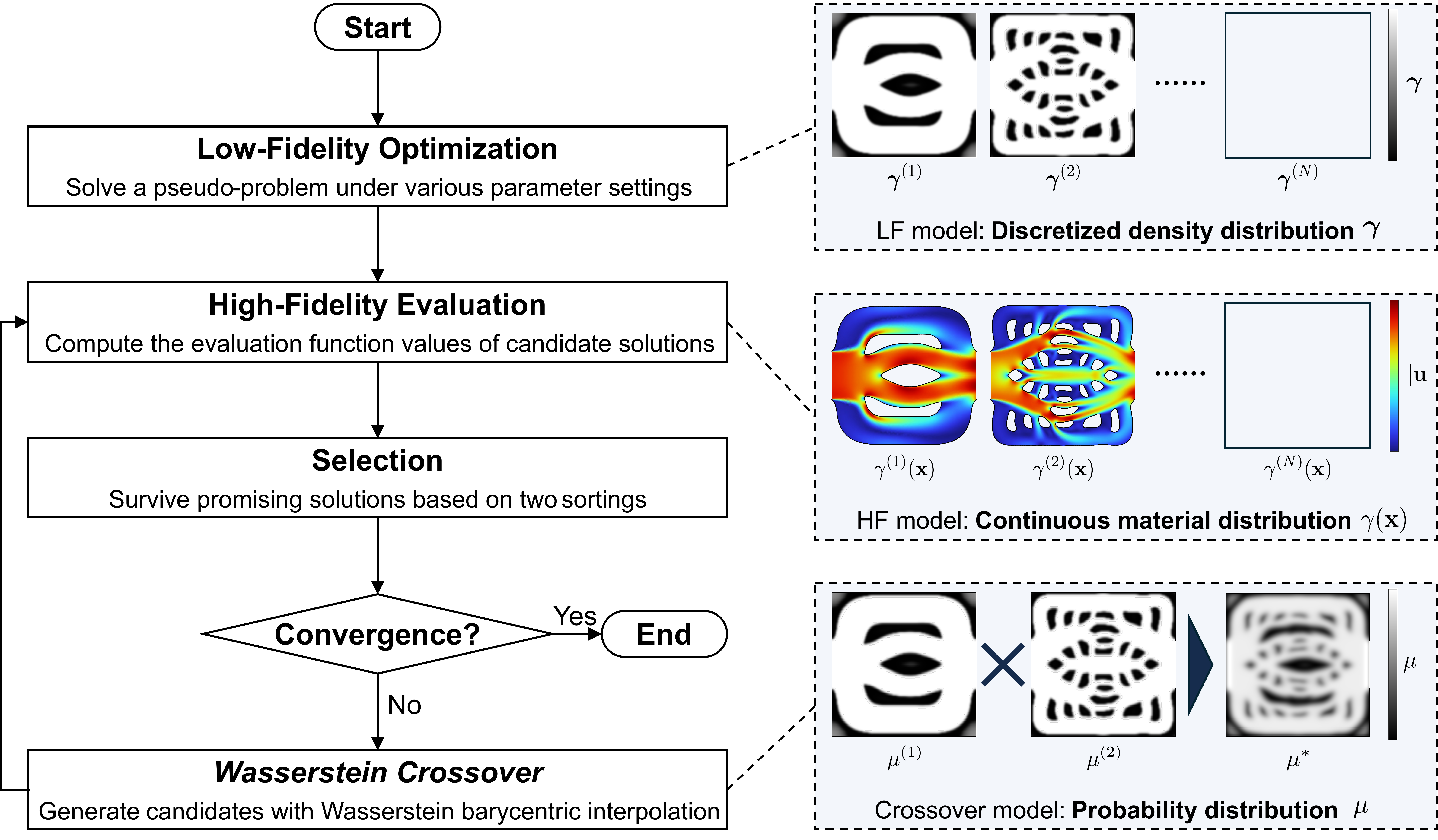} 
  \caption{Schematic flowchart of the proposed framework and representation of material distributions in each process}\label{fig_flowchart}
\end{figure}

Based on the above formulation, we describe each procedure of the proposed framework in the following.
The overall flowchart of the proposed framework is illustrated in Fig.~\ref{fig_flowchart}, whose main components are iterative evaluation, selection, and crossover.
The multifidelity formulations of Eq.~\eqref{eq_hf} and Eq.~\eqref{eq_lf} play a supplementary role in solving intractable topology optimization problems efficiently.
In the following, we describe the details of each procedure in the proposed framework, including LF optimization, HF evaluation, selection, convergence check, and Wasserstein crossover.
It should be noted that mutation, which is implemented with LF optimization in DDTD proposed by Yaji et al.~\cite{yaji2022}, is not incorporated into the proposed framework because its impact in DDTD has been observed to be limited, while it also incurs additional computational cost.
The development of effective mutation strategies is left as a subject for future work.

\subsection*{Low-fidelity optimization}\label{sec31}

The first step of the proposed framework is to solve the LF optimization problem of Eq.~\eqref{eq_lf}.
The primary aim of this step is to derive a diverse and promising set of density distributions, which serves as the initial population for the subsequent evolutionary process.
The key component of its role is the seeding parameter $\mathbf{s}$.
For example, by varying optimization algorithmic parameters such as filter radius or projection smoothness and threshold, the LF optimization problem can be solved multiple times to derive a diverse set of optimized density distributions.
In addition, employing different weightings in the weighted sum method or varying constraint thresholds in the $\varepsilon$-constraint method enables the derivation of certain results that reflect different trade-offs between the objectives, thereby providing a promising set of candidate solutions.

\subsection*{High-fidelity evaluation}\label{sec32}

In this step, each candidate solution is evaluated under the HF model to compute the objective and constraint values in the original optimization problem of Eq.~\eqref{eq_hf}.
A crucial aspect of this step is that it only requires the numerical evaluation of those functions, without the need for sensitivity analysis concerning the design variables.
This aspect allows the proposed framework to accommodate arbitrary evaluation functions, including non-differentiable and multimodal ones, thereby enabling the handling of minimax formulations and other complex optimization problems.
Moreover, this property allows each candidate represented as a discretized density distribution $\boldsymbol{\gamma}$ to be converted into a continuous material distribution $\gamma(\mathbf{x})$, remeshed, and then evaluated under the HF model.

\subsection*{Selection}\label{sec33}

In this step, promising candidate solutions are retained for the next generation, while the others are eliminated from the population based on their performance values obtained in the previous step.
Selection plays a crucial role in EAs, as it helps maintain population diversity and facilitates effective global search.
Since the proposed framework targets the multi-objective optimization problem defined in Eq.~\eqref{eq_hf}, ensuring sufficient diversity becomes even more challenging. 
To address this, we adopt the selection strategy in DDTD proposed by Kii et al.~\cite{kii2025}, which is designed to preserve the intrinsic diversity of material distributions.
This procedure is based on the non-dominated sorting genetic algorithm II (NSGA-II)~\cite{deb2002, verma2021}, which is one of the representative multi-objective GAs, and ranks candidates according to a two-stage sorting strategy.
The first stage is the non-dominated sorting, which assigns ranks to candidates based on their Pareto dominance relations.
The second stage, unique to their approach, further sorts candidates within the same rank using the Wasserstein distance sorting.
This is achieved by employing a topological data analysis method called persistent homology~\cite{edelsbrunner2002,zomorodian2005} to extract topological features from the material distributions, which are then quantified using the Wasserstein distance~\cite{mileyko2011}.
While the detailed explanation, including the implementation of these procedures, is left to the original paper~\cite{kii2025}, a predefined number of higher-ranked candidates are selected based on these two types of sorting, and the population for the current generation is formed accordingly.

\subsection*{Convergence check}\label{sec34}

In this step, the convergence of the optimization process is checked based on the hypervolume indicator~\cite{shang2020}, which is a representative performance metric for multi-objective optimization.
The hypervolume indicator quantifies the area or volume formed by the solution set $\mathcal{A}$ and a given reference point $\mathbf{r}\in\mathbb{R}^{n_{\text{o}}}$, defined as follows:
\begin{equation}\label{eq_hypervolume}
  \text{HV}(\mathcal{A}, \mathbf{r}) = \mathcal{L} \left(\bigcup _{\mathbf{a}\in\mathcal{A}}\{\mathbf{b}\,|\,\mathbf{a}\preceq\mathbf{b}\preceq\mathbf{r}\}\right),
\end{equation}
where $\mathcal{L}(\cdot )$ denotes the Lebesgue measure.
In this paper, we set the reference point $\mathbf{r}$ as 1.1 times the worst objective values among all the candidates in the initial population, which is a common practice in the literature~\cite{shang2020}.
The optimization is terminated either when the increase in hypervolume is deemed sufficiently converged or when a predefined number of iterations is reached.

\subsection*{Wasserstein crossover}\label{sec35}

This step constitutes the crossover operation proposed in this study for integration into DDTD, where new offspring are generated from parent individuals using interpolation based on the Wasserstein barycenter.
A key consideration here is how to interpret the discretized density distribution $\boldsymbol{\gamma}$ on the LF model in Eq.~\eqref{eq_lf} as a probability distribution to compute the Wasserstein barycenter in Eq.~\eqref{eq_barycenter}.
To this end, the density vector $\boldsymbol{\gamma}$ is first normalized into a probability vector $\mathbf{p} \in \mathbb{R}^n$ as follows:
\begin{equation}\label{eq_normalization}
  \mathbf{p} = \frac{\boldsymbol{\gamma}}{\sum_{e=1}^n \gamma_e}.
\end{equation}
Then, the corresponding discrete probability measure $\mu$ is defined by
\begin{equation}\label{eq_prob}
  \mu = \sum_{e=1}^n p_e\, \delta_{\mathbf{x}_e},
\end{equation}
where $\delta_{\mathbf{x}_e}$ denotes the Dirac measure at the $i$-th grid point $\mathbf{x}_e$, and $p_e$ represents the probability mass assigned to $\mathbf{x}_e$, i.e., $\mu(\{\mathbf{x}_e\}) = p_e$.
Based on this method of interpreting the discretized density distribution as a probability distribution, the proposed Wasserstein crossover takes the following steps:
\begin{enumerate}
  \item Randomly select two parents $\boldsymbol{\gamma}^{(i)}, \boldsymbol{\gamma}^{(j)} \in \varTheta^{(t)}$, where $\varTheta^{(t)}$ is the current population at iteration $t$.
  \item Convert the parents into probability distributions $\mu^{(i)}$ and $\mu^{(j)}$ using Eqs.~\eqref{eq_normalization} and \eqref{eq_prob}.
  \item Compute the Wasserstein barycenter $\mu^*$ of the two parental distributions $\mu^{(i)}$ and $\mu^{(j)}$ as follows:
  \begin{equation}\label{eq_barycenter2}
    \mu^* = \arg\min_{\mu} \left\{ \lambda W_{p, \varepsilon}(\mu, \mu^{(i)}) + (1-\lambda) W_{p, \varepsilon}(\mu, \mu^{(j)}) \right\},
  \end{equation}
  where $\lambda \in [0, 1]$ is a randomly sampled weight parameter.
  \item Convert the barycenter $\mu^* = \sum_{e=1}^n p_e^* \delta_{\mathbf{x}_e}$ back into a density distribution $\boldsymbol{\gamma}^*$ by applying min-max scaling to the probability vector $\mathbf{p}^* = [p_1^*, p_2^*, \dots, p_n^*]^\top$ as follows:
  \begin{equation}\label{eq_convert}
    \gamma_e^* = \frac{p_e^* - \min_{e'} p_{e'}^*}{\max_{e'} p_{e'}^* - \min_{e'} p_{e'}^*}, \quad \text{for } e = 1, 2, \dots, n,
  \end{equation}
  where $p_e^* = \mu^*(\{\mathbf{x}_e\})$ denotes the probability mass at the $e$-th grid point in the barycenter $\mu^*$.
\end{enumerate}
These sequential steps are repeated for a predefined number of offspring $N_{\text{xo}}$, each time adding the resulting density distribution $\boldsymbol{\gamma}^*$ as a new candidate solution to the population.

\section{Numerical implementation}\label{sec4}

In this section, we describe the detailed implementation of the proposed framework.
The overall procedure of the proposed framework, shown in Fig.~\ref{fig_flowchart}, is summarized in Algorithm~\ref{alg_framework}.
Note that the LF optimization and HF evaluation steps in Algorithm~\ref{alg_framework} are independent of the others, allowing for acceleration through parallel implementation as computational resources permit.
The entire algorithm was implemented using MATLAB (version 2023b).
In the following, we describe the implementation details of each operation in Algorithm~\ref{alg_framework}.

\begin{algorithm}[ht]
  \caption{Overall procedure of the proposed framework}\label{alg_framework}
  \begin{algorithmic}[1]
    \For{$k=1$ to $N_{\text{lf}}$}
      \State Solve the LF optimization problem of Eq.~\eqref{eq_lf} for $\boldsymbol{\gamma}^{(k)}$ on $\mathbf{s}^{(k)}$
    \EndFor
    \State Assemble a temporary dataset of the LF optimized density distributions, $\varTheta_{\text{tmp}}\gets\{\boldsymbol{\gamma}^{(1)}, \ldots, \boldsymbol{\gamma}^{(N_{\text{lf}})}\}$
    \For{$t=0$ to $t_{\text{max}}$}
      \For{$k=1$ to $|\varTheta_{\text{tmp}}|$}
        \State Calculate HF evaluation function values in Eq.~\eqref{eq_hf}, $\{J_1^{(k)}, \ldots, J_{n_{\text{o}}}^{(k)}\}$ and $\{G_1^{(k)}, \ldots, G_{n_{\text{c}}}^{(k)}\}$
      \EndFor
      \State Remove candidates from $\varTheta_{\text{tmp}}$ that violate the constraints $G_j^{(k)}\leq 0 \text{ for all } j=1, \ldots, n_{\text{c}}$
      \If{$t=0$}
        \State Set the current candidates, $\varTheta^{(t)}\gets\varTheta_{\text{tmp}}$
      \Else
        \State Augment the current candidates, $\varTheta^{(t)}\gets\varTheta^{(t-1)}\cup\varTheta_{\text{tmp}}$
      \EndIf
      \State Update $\varTheta^{(t)}$ using the selection algorithm based on $\mathcal{J}(\varTheta^{(t)})\gets\bigcup_{k=1}^{|\varTheta^{(t)}|}\{J_1^{(k)}, \ldots, J_{n_{\text{o}}}^{(k)}\}$
      \If{$\text{HV}$ in Eq.~\eqref{eq_hypervolume} regarding $\mathcal{J}(\varTheta^{(t)})$ is converged}
        \State \textbf{break}
      \EndIf
      \State Initialize a temporary dataset, $\varTheta_{\text{tmp}}\gets\emptyset$
      \For{$k=1$ to $N_{\text{xo}}$}
        \State Randomly select two parents $\boldsymbol{\gamma}^{(i)}, \boldsymbol{\gamma}^{(j)}\in\varTheta^{(t)}$
        \State Generate an offspring $\boldsymbol{\gamma}^{(k)}$ as the Wasserstein barycenter of $\boldsymbol{\gamma}^{(i)}$ and $\boldsymbol{\gamma}^{(j)}$
        \State Add the offspring to the temporary dataset, $\varTheta_{\text{tmp}}\gets\varTheta_{\text{tmp}}\cup\{\boldsymbol{\gamma}^{(k)}\}$
      \EndFor
    \EndFor
    \State \Return the optimized solutions $\varTheta^{(t)}$ and their performance values $\mathcal{J}(\varTheta^{(t)})$
  \end{algorithmic}
\end{algorithm}

\subsection{Low-fidelity optimization}\label{sec41}

The LF optimization problem of Eq.~\eqref{eq_lf} is solved using the density-based method~\cite{bendsoe1989}, which is a representative approach for topology optimization.
To ensure the smoothness of the design variable $\boldsymbol{\gamma}^{(k)}$ in Eq.~\eqref{eq_lf}, we employ a density filter~\cite{bruns2001,bourdin2001} defined as follows:
\begin{equation}\label{eq_filter}
  \widetilde{\gamma}_e^{(k)} = \frac{\sum_{e'\in\mathcal{N}_e} w_{ee'} \gamma_{e'}^{(k)}}{\sum_{e'\in\mathcal{N}_e} w_{ee'}}, \quad \text{for } e=1, 2, \ldots, n,
\end{equation}
where $\mathcal{N}_e$ is the set of neighboring elements of the $e$-th element defined by:
\begin{equation}\label{eq_ele}
  \mathcal{N}_e = \{e'\,|\,\|\mathbf{x}_e - \mathbf{x}_{e'}\|_2 \leq R\},
\end{equation}
and $w_{ee'}$ is the weight assigned to element $e'$ with respect to element $e$, given by:
\begin{equation}\label{eq_weight}
  w_{ee'} = 1 - \frac{\|\mathbf{x}_e - \mathbf{x}_{e'}\|_2}{R},
\end{equation}
where $R$ is the filter radius.
We use the method of moving asymptotes (MMA)~\cite{svanberg1987} for updating the design variables, and the move limit is set to 0.05 for all calculations of LF optimization in this paper.
The entire LF optimization process was implemented in MATLAB, while the finite element analysis and sensitivity analysis, based on the discrete adjoint method, were performed using COMSOL Multiphysics (version 6.3).

\subsection{High-fidelity evaluation}\label{sec42}

\begin{figure}[t]
  \centering
  \includegraphics[width=.9\textwidth]{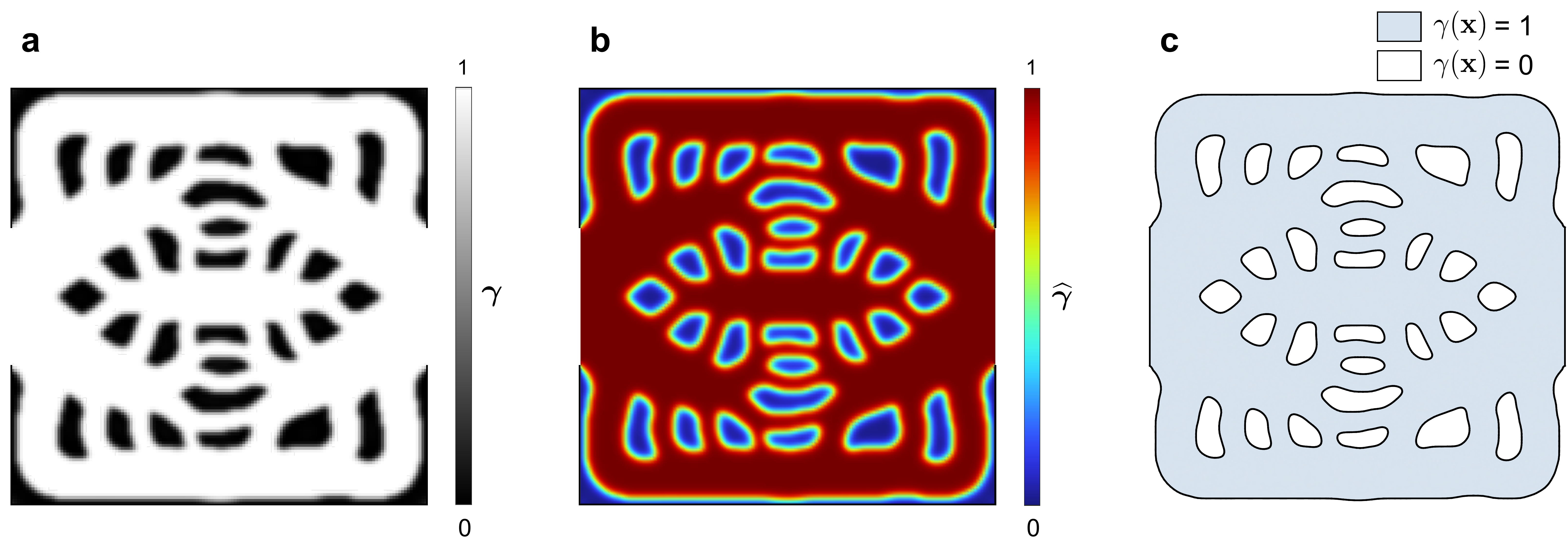}
  \caption{Pre-process for high-fidelity evaluation: (a) the density distribution $\boldsymbol{\gamma}$ on the low-fidelity model; (b) the filtered density distribution $\widehat{\boldsymbol{\gamma}}$; (c) the material distribution $\gamma(\mathbf{x})$ on the high-fidelity model.}\label{hfproc}
\end{figure}

To compute the evaluation function values in the HF model of Eq.~\eqref{eq_hf}, it is necessary to convert the discretized density distribution $\boldsymbol{\gamma}$ on the LF model into a continuous material distribution $\gamma(\mathbf{x})$ on the HF model as a pre-processing step.
The pre-processing procedure adopted in this paper is illustrated in Fig.~\ref{hfproc}.
First, we compute the filtered density distribution $\widehat{\boldsymbol{\gamma}}$ using a partial differential equation (PDE)-based filter~\cite{lazarov2011} defined as follows:
\begin{equation}\label{eq_pdefil}
  -R_{\text{h}}^2\nabla^2 \widehat{\boldsymbol{\gamma}} + \widehat{\boldsymbol{\gamma}} = \boldsymbol{\gamma},
\end{equation}
where $R_{\text{h}}$ is the filter radius.
Note that when solving the PDE in Eq.~\eqref{eq_pdefil}, Dirichlet boundary conditions are imposed on specific regions as necessary: for example, fixed and loaded boundaries in structural mechanics problems, or inlet and outlet boundaries in thermofluid problems, where  $\widehat{\boldsymbol{\gamma}}_e = 1$ or $\widehat{\boldsymbol{\gamma}}_e = 0$ is prescribed accordingly.
An isosurface at $\widehat{\boldsymbol{\gamma}}_e = 0.5$ is then extracted from the filtered distribution field.
This process yields a smooth boundary, from which a continuous material distribution $\gamma(\mathbf{x})$ is obtained.
The resulting material distributions are then discretized using appropriate meshes based on the problem settings, and the HF evaluation function values are subsequently evaluated.
All of these HF evaluation processes were implemented in COMSOL Multiphysics.

\subsection{Wasserstein crossover}\label{sec43}

In this paper, we use the 2-Wasserstein distance with entropic regularization in Eq.~\eqref{eq_barycenter2} for the barycentric computation, and replace certain matrix operations in Algorithm~\ref{alg_barycenter} with convolutions using a Gaussian filter to accelerate the computation as described in Section~\ref{sec22}.
In addition, the entropic regularization coefficient $\varepsilon$, which controls the trade-off between the accuracy of optimal transport and computational cost, is adaptively adjusted based on the pair of parent distributions selected for crossover.
Specifically, prior to the crossover step described in Section~\ref{sec35}, we compute a distance matrix $\mathbf{D}\in\mathbb{R}^{N_{\text{pop}}\times N_{\text{pop}}}$, where $N_{\text{pop}}$ is the population size, defined as follows:
\begin{equation}\label{eq_dist}
  D_{ij}=\|\boldsymbol{\gamma}^{(i)} - \boldsymbol{\gamma}^{(j)}\|_2 = \sqrt{\sum_{e=1}^n \left(\gamma_e^{(i)} - \gamma_e^{(j)}\right)^2}.
\end{equation}
Using the maximum and minimum value $D_{\max}, D_{\min}$ of the distance matrix $\mathbf{D}$, the entropic regularization coefficient $\varepsilon$ for crossover between two parents $\boldsymbol{\gamma}^{(i)}$ and $\boldsymbol{\gamma}^{(j)}$ is adaptively determined as follows:
\begin{equation}\label{eq_adaptive}
  \varepsilon = \varepsilon_{\min} + (\varepsilon_{\max} - \varepsilon_{\min}) \frac{D_{ij} - D_{\min}}{D_{\max} - D_{\min}},
\end{equation}
where $\varepsilon_{\max}$ and $\varepsilon_{\min}$ are the maximum and minimum values of the entropic regularization coefficient, respectively.
In this paper, we determine $\varepsilon_{\max}$ and $\varepsilon_{\min}$, and the convergence tolerance $\tau$ empirically based on preliminary numerical experiments for each problem setting, taking into account the trade-off between the degree of blurring introduced by regularization, the validity of the interpolated results, and the computational cost.
By adaptively adjusting the entropic regularization coefficient in this manner, the Wasserstein barycenter can be computed more accurately when the parents are relatively similar within the population, while maintaining computational efficiency when they are dissimilar, which enhances the crossover performance.
The Wasserstein crossover procedure was entirely implemented in Python (version 3.10.16).
To accelerate the computation of the algorithm described in Algorithm~\ref{alg_barycenter}, we employed JAX (version 0.6.1)~\cite{bradbury2018} and implemented it based on POT (version 0.9.5)~\cite{flamary2021}, a Python library for numerical computation of optimal transport.

\section{Numerical examples}\label{sec5}

In this section, we demonstrate the effectiveness of the proposed framework for intractable topology optimization problems through three numerical examples.
The first example is a two-dimensional stress minimization problem that includes a non-differentiable objective function in the formulation.
We verify the effectiveness of the framework for structural mechanics problems and highlight the advantages of the proposed Wasserstein crossover by comparing it with the VAE-based crossover employed in conventional DDTD.
The second example is a two-dimensional turbulent heat transfer problem, where evaluation functions exhibit strong multimodality due to the complex physics of turbulent flow.
This example illustrates the applicability of the proposed framework beyond simple structural mechanics problems, extending it to complex thermofluid problems.
We tackle a three-dimensional stress minimization problem as the final example, demonstrating the suitability for three-dimensional topology optimization problems.
Finally, we discuss the computational cost for each example to evaluate the scalability of the proposed framework.

All numerical examples presented in this paper were computed on a Linux workstation with a 2.7 GHz AMD Ryzen Threadripper PRO 3995WX 64-core processor with 512 GB RAM and an NVIDIA RTX A6000 GPU with 48 GB memory.
For the two-dimensional cases, the LF optimization and HF evaluation processes were executed in parallel using 50 CPU cores, whereas 30 CPU cores were used for the three-dimensional case.
In all examples, the Wasserstein crossover process was computed on the GPU.

\begin{figure}[t]
  \centering
  \includegraphics[width=\textwidth]{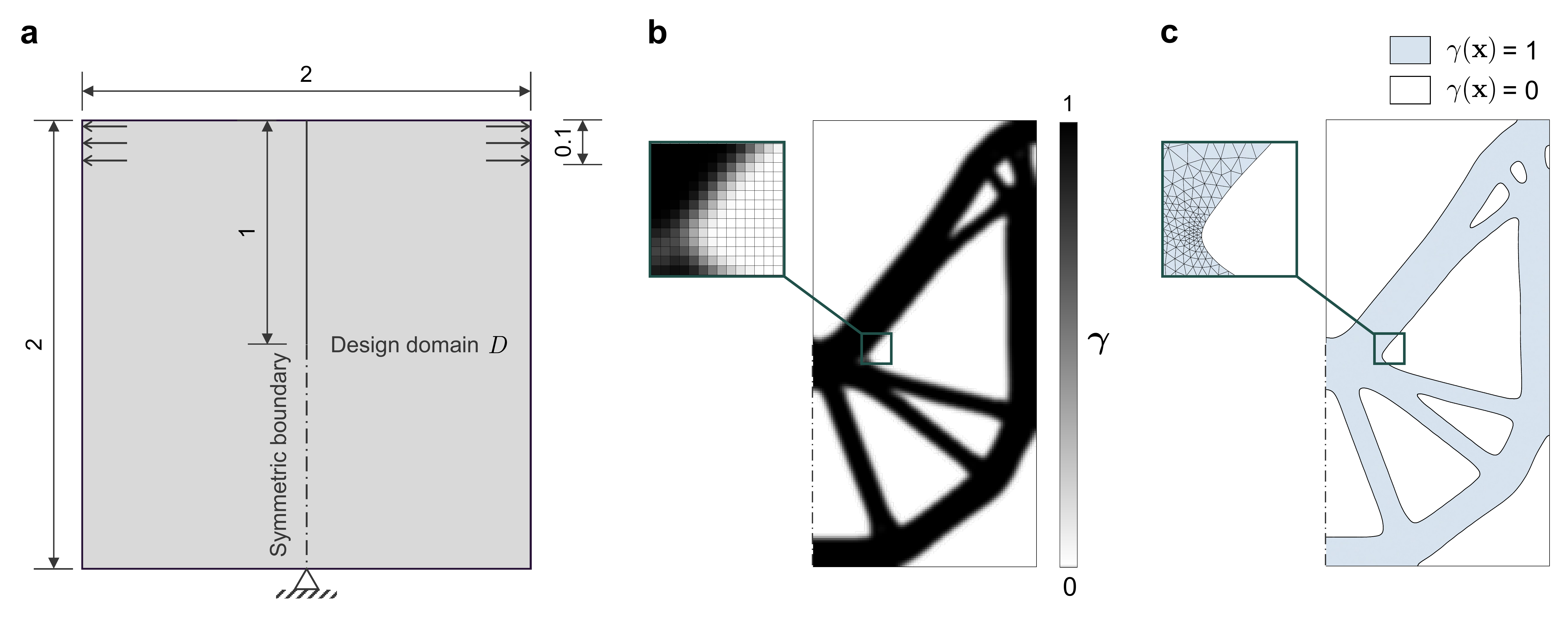}
  \caption{Problem settings in the two-dimensional cracked plate design problem: (a) the design domain and boundary conditions; (b, c) a topology optimized configuration with (b) a density distribution $\boldsymbol{\gamma}$ discretized with structured meshes on the low-fidelity model, and (c) a material distribution $\gamma(\mathbf{x})$ discretized with body-fitted meshes on the high-fidelity model.}\label{problem_crack}
\end{figure}

\subsection{2D stress minimization problem}\label{sec51}

\subsubsection{Problem settings}\label{sec511}

Let us consider a multi-objective topology optimization problem of a two-dimensional cracked plate design, which is often used as a benchmark for stress-based topology optimization~\cite{emmendoerfer2014, emmendoerfer2016, chu2018, giraldo-londono2021, norato2022}.
The optimization problem is formulated as follows:
\begin{eqnarray}\label{eq_crackhf}
  \begin{aligned}
    & \underset{\gamma(\mathbf{x})}{\text{minimize}}
    && J_1 = \underset{\mathbf{x}\in D}{\text{max}} \left\{ \sigma(\mathbf{x}) \right\}, \\
    &&& J_2 = \frac{\int_D \gamma(\mathbf{x})\,\mathrm{d}\mathit{\Omega}}{\int_D \mathrm{d}\mathit{\Omega}}, \\
    & \text{subject to}
    && \gamma(\mathbf{x})\in\{0, 1\}, \quad \forall\mathbf{x}\in D.
  \end{aligned}
\end{eqnarray}
Herein, $\sigma(\mathbf{x})$ is the von Mises stress at the point $\mathbf{x}\in D \subset \mathbb{R}^2$.
It is formulated as a two-objective optimization problem of minimizing the maximum stress $J_1$ and volume fraction $J_2$.
While the formulation of Eq.~\eqref{eq_crackhf} is considered as the HF model corresponding to the formulation of Eq.~\eqref{eq_hf}, the LF optimization problem corresponding to Eq.~\eqref{eq_lf} is formulated as follows:
\begin{eqnarray}\label{eq_cracklf}
  \begin{aligned}
    & \underset{\boldsymbol{\gamma}^{(k)}}{\text{minimize}}
    && \widetilde{J} = \left(\sum_{e=1}^n \sigma_e^P\right)^{1/P}, \\
    & \text{subject to}
    && \widetilde{G} = \sum_{e=1}^n v_e\gamma_e^{(k)} - V\sum_{e=1}^n v_e\leq 0, \\
    &&& \boldsymbol{\gamma}^{(k)}\in [0, 1]^n, \\
    & \text{for given } && \mathbf{s}^{(k)}=[s_1, s_2]^\top.
  \end{aligned}
\end{eqnarray}
Herein, $\sigma_e$ and $v_e$ are the elemental von Mises stress and volume, respectively.
Two seeding parameters, $s_1$ and $s_2$, are used to determine the filter radius $R$ in Eq.~\eqref{eq_filter} and volume constraint $V$, respectively.
Each parameter is sampled uniformly from the interval $[0,1]$ with $n_{s_1}$ and $n_{s_2}$ divisions.
The filter radius $R$ is obtained by linearly interpolating between $R_{\text{min}}$ and $R_{\text{max}}$ according to $s_1$, while the volume constraint $V$ is similarly interpolated between $V_{\text{min}}$ and $V_{\text{max}}$ using $s_2$, resulting in a total of $N_{\text{lf}} = n_{s_1} \times n_{s_2}$ parameter combinations.
The objective function $\widetilde{J}$ in the LF formulation of Eq.~\eqref{eq_cracklf} is defined as the $P$-norm stress, which is a differentiable approximation of the maximum stress.
For the norm parameter $P$, we set a constant value of $P=8$, which is often used from the perspective of the trade-off between accuracy and numerical stability~\cite{le2010}.
The LF optimization problem of Eq.~\eqref{eq_cracklf} is solved independently $N_{\text{lf}}$ times using the solid isotropic material with penalty (SIMP) method~\cite{bendsoe2003}.

\begin{table}[t]
  \centering
  \caption{Overall parameters in the two-dimensional stress minimization problem.}\label{tab_crack}
  \begin{tabular*}{\textwidth}{@{\extracolsep\fill}lll}
    \toprule
    Description & Symbol & Value \\
    \midrule
    Maximum number of iterations & $t_{\text{max}}$ & 100 \\
    Population size & $N_{\text{pop}}$ & 100 \\ 
    Number of offspring & $N_{\text{xo}}$ & 100 \\
    Number of LF optimized initial designs & $N_{\text{lf}}$ & 100 \\
    Number of seeding parameters for filter radius & $n_{s_1}$ & 4 \\
    Number of seeding parameters for volume fraction & $n_{s_2}$ & 25 \\
    Minimum filter radius & $R_{\text{min}}$ & $0.03$ \\
    Maximum filter radius & $R_{\text{max}}$ & $0.12$ \\
    Minimum volume fraction & $V_{\text{min}}$ & $0.30$ \\
    Maximum volume fraction & $V_{\text{max}}$ & $0.60$ \\
    Filter radius for HF evaluation & $R_{\text{h}}$ & $0.01$ \\
    Minimum entropic regularization coefficient & $\varepsilon_{\min}$ & $1.0\times 10^{-6}$ \\
    Maximum entropic regularization coefficient & $\varepsilon_{\max}$ & $1.0\times 10^{-4}$ \\
    Convergence tolerance in Sinkhorn algorithm & $\tau$ & $1.0\times 10^{-9}$ \\
    \bottomrule
  \end{tabular*}
\end{table}

\begin{table}[t]
  \centering
  \caption{Parameters for the VAE-based crossover in the two-dimensional stress minimization problem.}\label{tab_crackvae}
  \begin{tabular*}{\textwidth}{@{\extracolsep\fill}lll}
    \toprule
    Description & Value \\
    \midrule
    Size of input and output data & 20,000 \\
    Size of latent space & 8 \\
    Architecture of encoder & [20,000, 512, 8] (fully connected) \\
    Architecture of decoder & [8, 512, 20,000] (fully connected) \\
    Activation function & Sigmoid (output layer)\\
    & ReLU (other layers) \\
    Reconstruction loss function & Mean squared error \\
    Weighting coefficient for Kullback-Leibler (KL) divergence & $0.01$ \\
    Optimizer & Adam~\cite{kingma2014} \\
    Learning rate & $0.001$ \\
    Batch size & 10 \\
    Number of epochs & 500 \\
    \bottomrule
  \end{tabular*}
\end{table}

The design domain and boundary conditions of the cracked plate are illustrated in Fig.~\ref{problem_crack}(a).
Due to the symmetry, we define the right half of a $2\times 2$ plate as the design domain $D$, where a symmetric boundary condition is imposed only on the lower half of the central boundary to represent a crack in the upper half.
In the LF model shown in Fig.~\ref{problem_crack}(b), the design variable field $\boldsymbol{\gamma}^{(k)}$ is discretized using structured meshes with $200\times 100$ quadrilateral elements, i.e., $n=20000$.
In contrast, body-fitted meshes are employed in the HF model shown in Fig.~\ref{problem_crack}(c), where the boundary between material and void regions is defined by an isosurface extracted through the pre-processing procedure described in Section~\ref{sec42}.
The maximum and minimum element sizes for the HF model are set to $0.04$ and $1.5\times 10^{-4}$, respectively, and body-fitted meshes with quadrilateral elements are automatically built for each candidate solution.
In the LF model, the design variable $\boldsymbol{\gamma}^{(k)}$ is discretized using $\mathbb{P}^0$ Lagrange finite elements, whereas the HF model employs $\mathbb{P}^2$ Lagrange finite elements for structural analysis.

The parameters regarding overall procedures of the proposed framework are listed in Table~\ref{tab_crack}.
Note that the convergence check with the hyper volume indicator is not performed in this example, and the optimization is terminated after a predefined number of iterations $t_{\text{max}}$.

Using this numerical example, we compare the VAE-based crossover employed in conventional DDTD and the proposed Wasserstein crossover.
The settings for the VAE in this comparison are listed in Table~\ref{tab_crackvae}. 
We determined these settings following the previous studies~\cite{kii2024, kato2025, kii2025} where DDTD was applied to stress-based topology optimization problems similar to this example.

\subsubsection{Results and discussion}\label{sec512}

\begin{figure}[t]
  \centering
  \includegraphics[width=\textwidth]{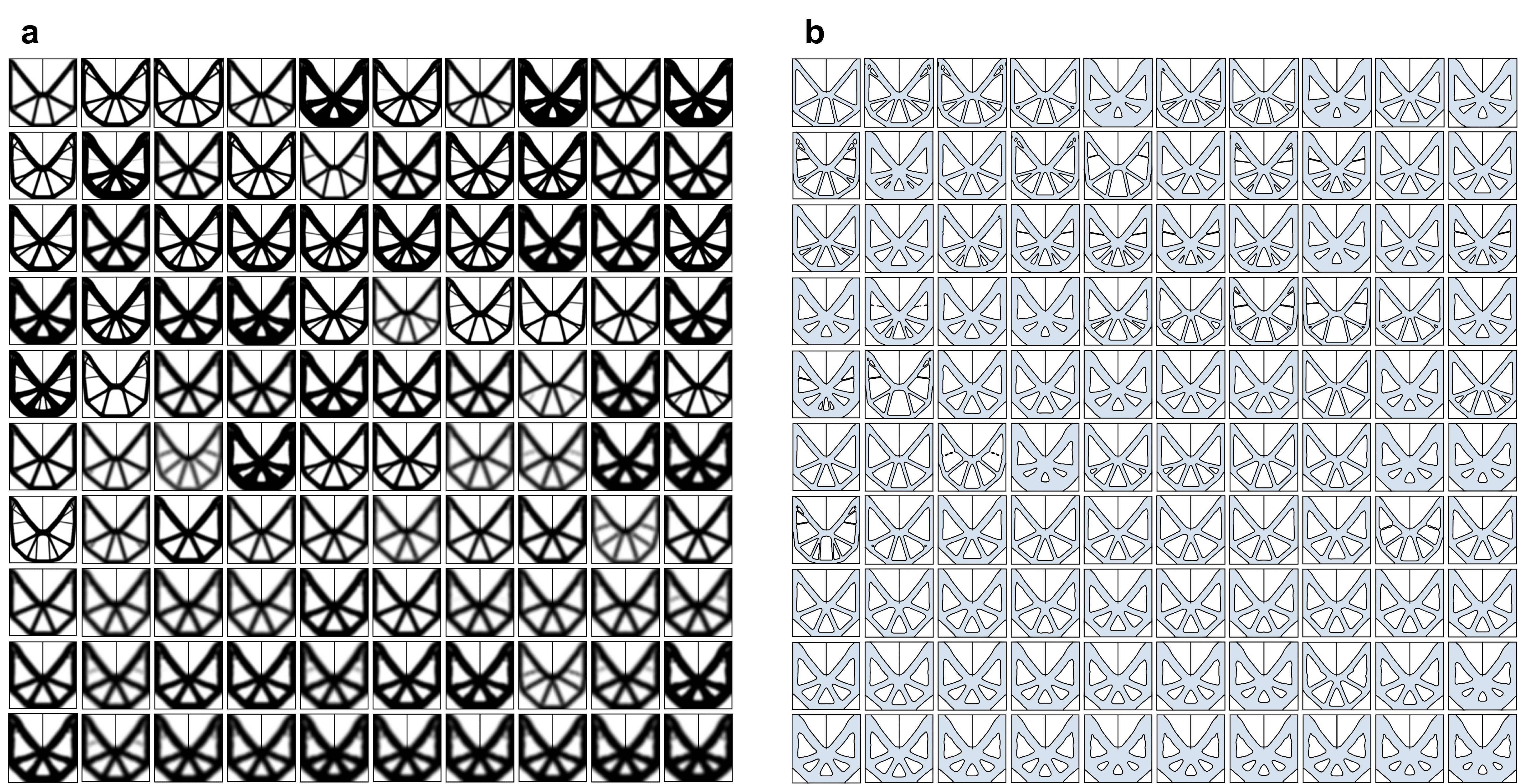}
  \caption{Low-fidelity optimized initial designs in the two-dimensional cracked plate design problem: (a) density distributions on the low-fidelity model; (b) material distributions on the high-fidelity model.}\label{result_crack_ini}
\end{figure}

Figure~\ref{result_crack_ini} shows the initial population obtained from the LF optimization problem of Eq.~\eqref{eq_cracklf}.
Owing to the effect of the seeding parameters $s_1$ and $s_2$ in Eq.~\eqref{eq_cracklf}, density distributions on the LF model shown in Fig.~\ref{result_crack_ini}(a) exhibit a variety of topologies under different volume fractions.
These density distributions are then converted into binarized material distributions on the HF model, as shown in Fig.~\ref{result_crack_ini}(b).
Some disconnected structures are observed in Fig.~\ref{result_crack_ini}(b), which results from the fact that the Heaviside projection technique~\cite{wang2011}, commonly used to promote binarization in the density-based method, was deliberately not applied to prioritize stability in the LF optimization.

\begin{figure}[t]
  \centering
  \includegraphics[width=.95\textwidth]{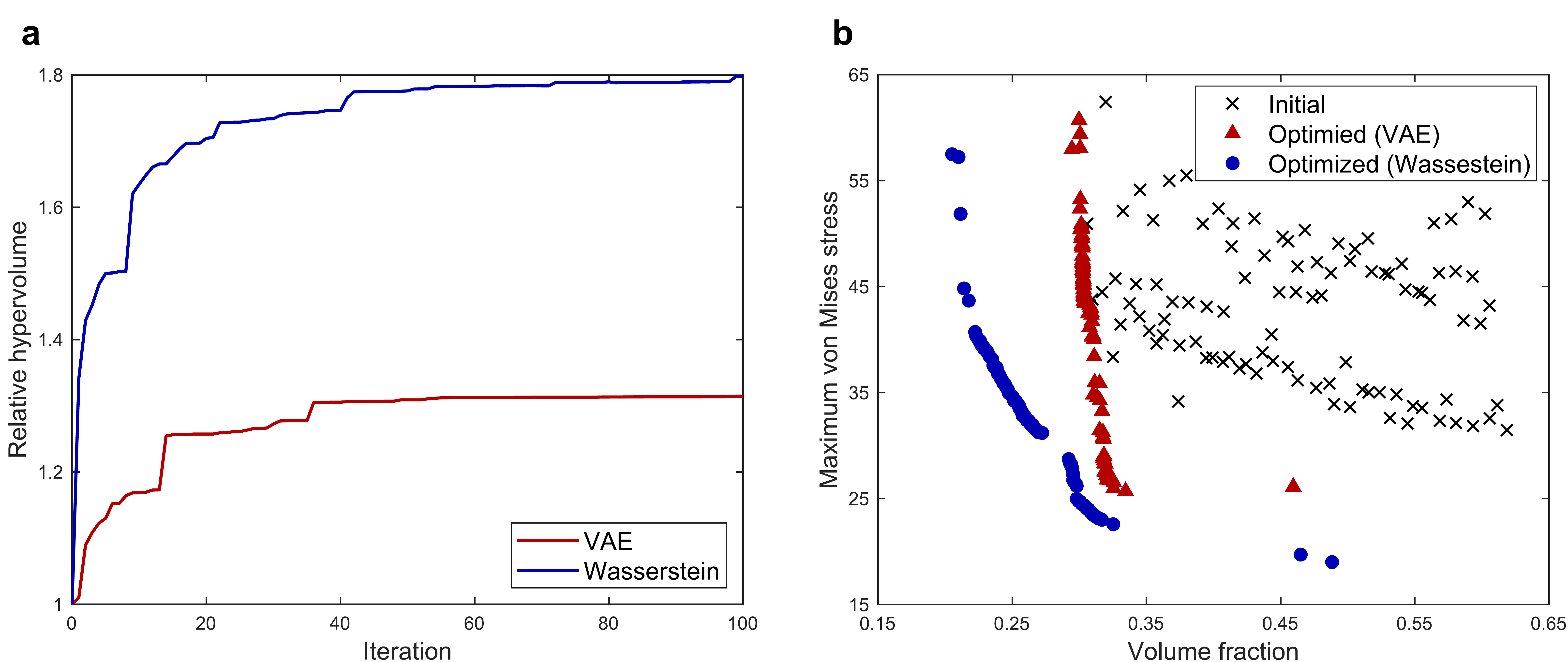}
  \caption{Optimization results in the two-dimensional cracked plate design problem: (a) convergence history of the hypervolume indicator normalized by the initial value; (b) objective space where the performance values of the initial and optimized designs by VAE-based and Wasserstein crossover are plotted.}\label{result_crack}
\end{figure}

\begin{figure}[t]
  \centering
  \includegraphics[width=\textwidth]{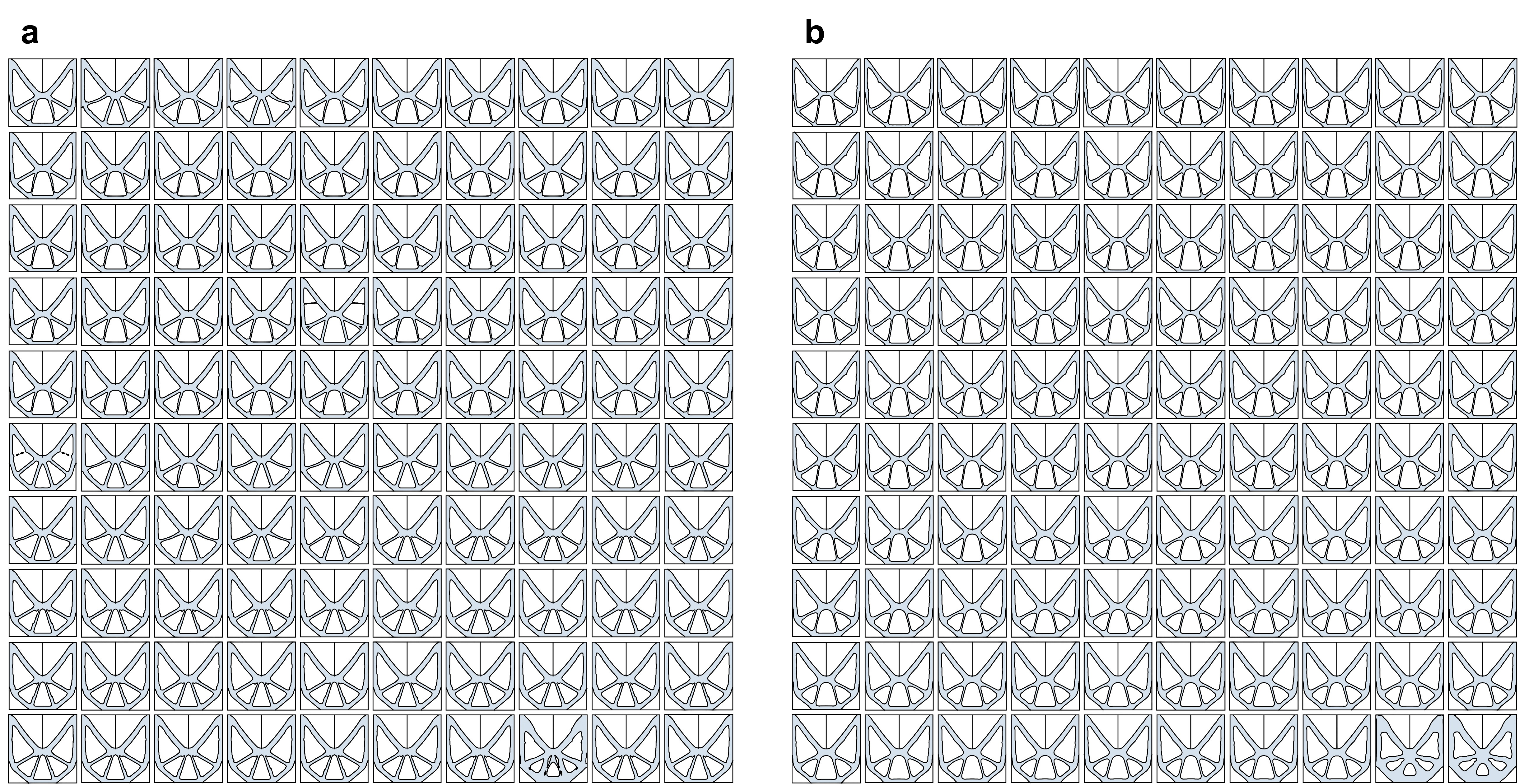}
  \caption{Optimized designs in the two-dimensional cracked plate design problem: (a) optimized via VAE-based crossover; (b) optimized via Wasserstein crossover.}\label{result_crack_opt}
\end{figure}

\begin{figure}[t]
  \centering
  \includegraphics[width=.77\textwidth]{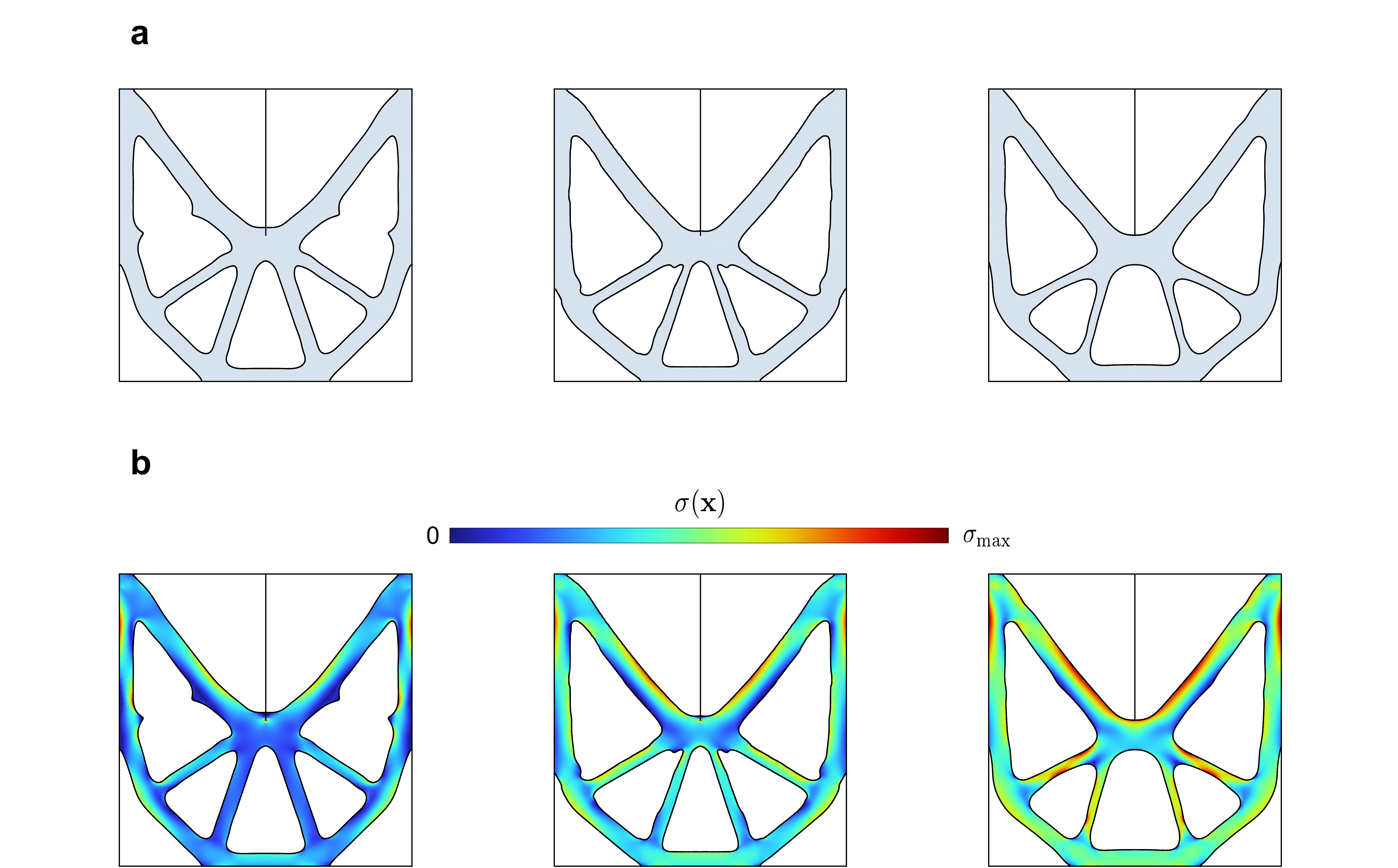}
  \caption{Results comparison of initial (left), optimized via VAE-based crossover (middle), and optimized via Wasserstein crossover (right) designs under nearly identical volume fraction conditions in the two-dimensional cracked plate design problem: (a) material distributions; (b) stress distributions $\sigma(\mathbf{x})$. Herein, these objective values are: $J_1=38.4$, $J_2=0.325$ (left); $J_1=27.5$, $J_2=0.320$ (middle); $J_1=23.0$, $J_2=0.317$ (right), where $J_1$ and $J_2$ are the maximum stress and volume fraction, respectively. Note that the color bar is scaled based on the maximum stress $\sigma_{\text{max}}$ of each design.}\label{result_crack_compare1}
\end{figure}

Figure~\ref{result_crack} shows the optimization population comparing the VAE-based crossover and the proposed Wasserstein crossover.
The initial designs shown in Fig.~\ref{result_crack_ini} are progressively evolved through iterations using the VAE-based and Wasserstein crossovers.
As shown in Fig.~\ref{result_crack}(a), Wasserstein crossover yields a significantly greater performance improvement.
Specifically, the VAE-based crossover achieves an approximate 31\% improvement from the initial value, whereas the Wasserstein crossover achieves an about 80\% improvement.
Figure~\ref{result_crack}(b) further illustrates that the Pareto front obtained via the Wasserstein crossover exhibits a substantial advancement from the initial solutions and completely dominates the Pareto front via the VAE-based crossover.

Figure~\ref{result_crack_opt} compares the optimized designs obtained via the VAE-based and Wasserstein crossovers.
A common feature observed in both cases is that the population has converged to material distributions with almost the same topology, which also exists in the initial population in Fig.~\ref{result_crack_ini}(b).
As can be seen in Fig.~\ref{result_crack_opt}(a), the optimized designs via Wasserstein crossover achieve a more significant reduction in volume.

\begin{figure}[t]
  \centering
  \includegraphics[width=.77\textwidth]{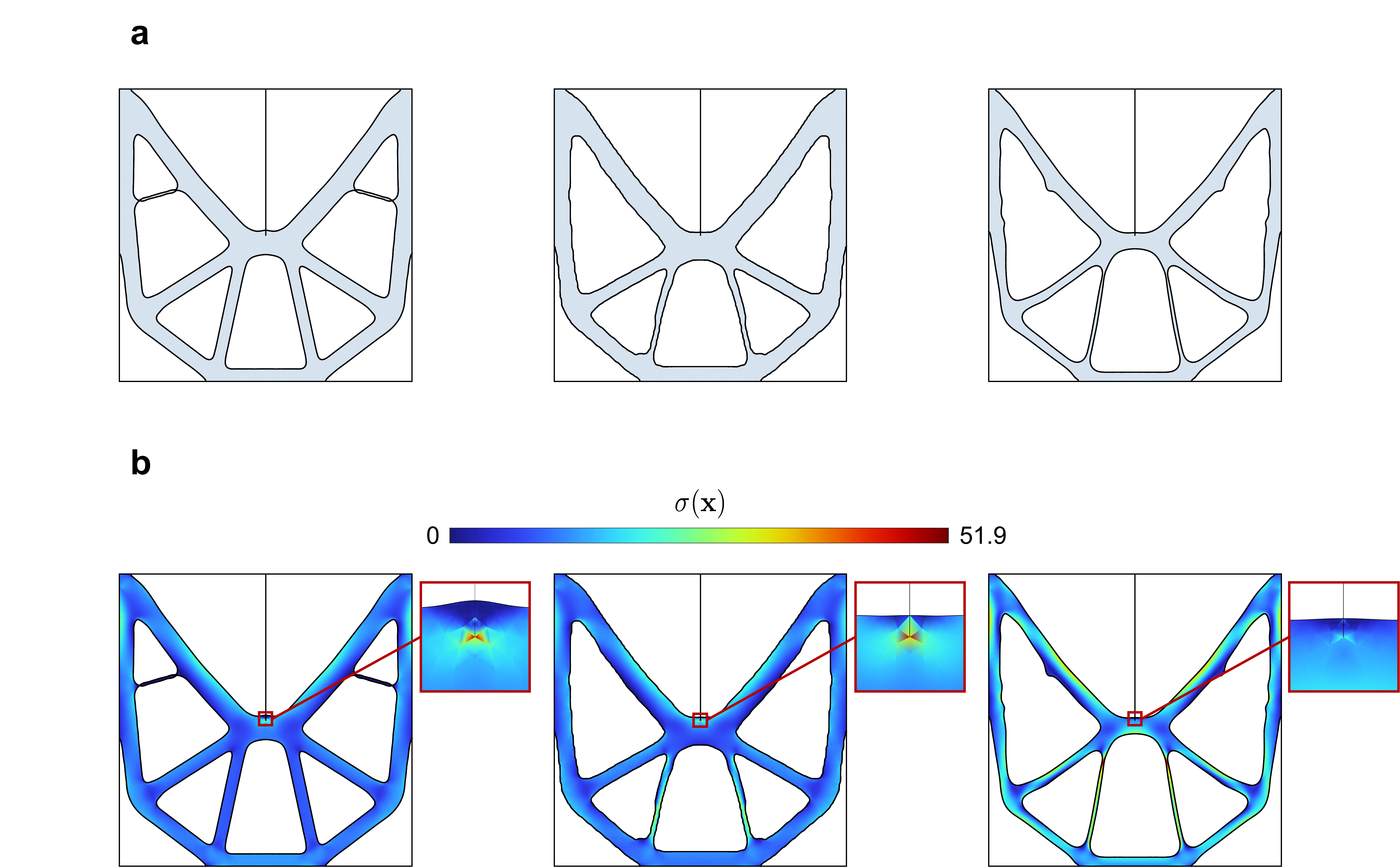}
  \caption{Results comparison of initial (left), optimized via VAE-based crossover (middle), and optimized via Wasserstein crossover (right) designs under nearly identical maximum stress conditions in the two-dimensional cracked plate design problem: (a) material distributions; (b) stress distributions $\sigma(\mathbf{x})$. Herein, these objective values are: $J_1=50.9$, $J_2=0.306$ (left); $J_1=50.9$, $J_2=0.301$ (middle); $J_1=51.9$, $J_2=0.212$ (right), where $J_1$ and $J_2$ are the maximum stress and volume fraction, respectively.}\label{result_crack_compare2}
\end{figure}

To enable a more detailed comparison, Fig.~\ref{result_crack_compare1} presents the comparison results of the initial, optimized via VAE-based crossover, and optimized via Wasserstein crossover designs under almost identical volume fraction conditions.
While these designs indeed exhibit the same topology, a comparison of the maximum stress value $J_1$ reveals that the optimized design via VAE-based crossover achieves an approximate 28\% reduction from the initial design, whereas the optimized design via Wasserstein crossover further improves upon this, achieving a 40\% reduction.
This improvement can be attributed to subtle geometric modifications, such as changes in the outer shape and the rounding of holes, which lead to more uniform stress distributions.
Notably, the optimized design via Wasserstein crossover particularly exhibits a more uniform stress distribution around the tip of the crack, which is typically prone to a high stress concentration.
Figure~\ref{result_crack_compare2} next presents a comparison under approximately identical conditions in terms of maximum stress $J_1$.
Although all three designs exhibit similar outer shapes, it can be seen that the thin structure, which is barely continuous in the initial design, is removed in both the optimized designs.
Focusing on the crack tip in Fig.~\ref{result_crack_compare2}(b), stress concentration is evident in the initial design and optimized design via VAE-based crossover.
In contrast, the optimized design via Wasserstein crossover effectively redistributes the stress throughout the surrounding thin structures, resulting in a reduction in the volume fraction $J_2$ by approximately 30\% compared to the initial design.

These results suggest that the proposed framework can be effectively applied to stress-based topology optimization problems, particularly in handling non-differentiable objective functions such as the maximum stress.
Furthermore, the proposed Wasserstein crossover appears capable of generating physically reasonable material distributions not only in fluid problems, as shown in Fig.~\ref{fig_morphing_compare}, but also in structural mechanics problems.
Its performance as a crossover operator seems to be superior to that of the VAE-based crossover, which is used in conventional DDTD.

\begin{figure}[t]
  \centering
  \includegraphics[width=\textwidth]{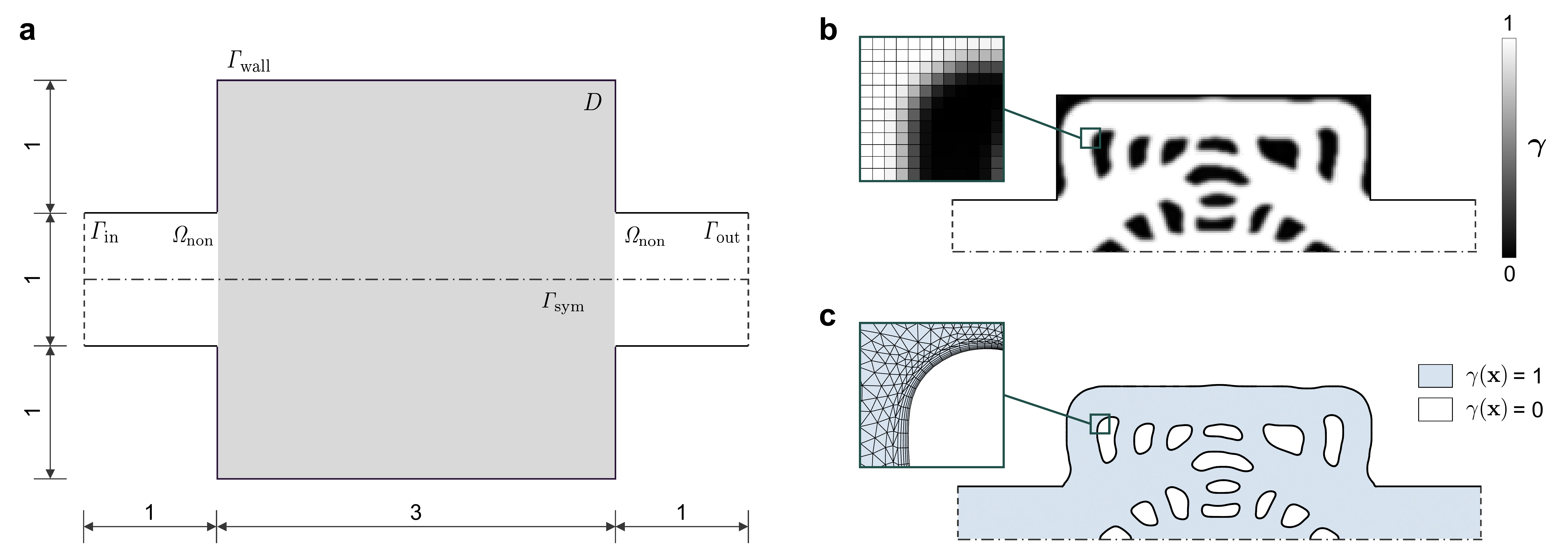}
  \caption{Problem settings in the two-dimensional heat sink design problem: (a) the analysis domain composed of the design domain $D$ and non-design domain $\mathit{\Omega}_{\text{non}}$; (b, c) a topology optimized configuration with (b) a density distribution $\boldsymbol{\gamma}$ discretized with structured meshes on the low-fidelity model, and (c) a material distribution $\gamma(\mathbf{x})$ discretized with body-fitted meshes and boundary layer meshes on the high-fidelity model. The fluid comes from the inlet $\mathit{\Gamma}_{\text{in}}$ to the outlet $\mathit{\Gamma}_{\text{out}}$. The boundary $\mathit{\Gamma}_{\text{sym}}$ is defined as a symmetric boundary. The remaining external boundary $\mathit{\Gamma}_{\text{wall}}$ is defined as a wall boundary.}\label{problem_thermofluid}
\end{figure}

\subsection{2D turbulent heat transfer problem}\label{sec52}

\subsubsection{Problem settings}\label{sec521}

Let us consider a heat sink design problem considering two-dimensional turbulent heat transfer.
Figure~\ref{problem_thermofluid}(a) illustrates the design domain $D\subset\mathbb{R}^2$ and an analysis domain $\mathcal{O}=D\cup\mathit{\Omega}_{\text{non}}$ with $D\cap \mathit{\Omega}_{\text{non}}\neq \emptyset$, where $\mathit{\Omega}_{\text{non}}$ is a non-design domain.
Examples of the optimized configuration on the LF and HF models are shown in Figs.~\ref{problem_thermofluid}(b) and (c), respectively.
The parameters regarding the overall algorithm are listed in Table~\ref{tab_thermofluid}, and the problem settings, including governing equations and formulations for both the HF and LF models, are described as follows:

\begin{table}[t]
  \centering
  \caption{Parameters regarding the overall algorithm in the two-dimensional turbulent heat transfer problem.}\label{tab_thermofluid}
  \begin{tabular*}{\textwidth}{@{\extracolsep\fill}lll}
    \toprule
    Description & Symbol & Value \\
    \midrule
    Maximum number of iterations & $t_{\text{max}}$ & 100 \\
    Population size & $N_{\text{pop}}$ & 100 \\ 
    Number of offspring & $N_{\text{xo}}$ & 100 \\
    Number of LF optimized initial designs & $N_{\text{lf}}$ & 100 \\
    Minimum entropic regularization coefficient & $\varepsilon_{\min}$ & $2.0\times 10^{-6}$ \\
    Maximum entropic regularization coefficient & $\varepsilon_{\max}$ & $2.0\times 10^{-4}$ \\
    Convergence tolerance in Sinkhorn algorithm & $\tau$ & $1.0\times 10^{-9}$ \\
    \bottomrule
  \end{tabular*}
\end{table}

\subsubsection*{High-fidelity model}
The analysis domain $\mathcal{O}$ on the HF model is composed of the fluid domain $\mathit{\Omega}_{\text{f}}$ and the solid domain $\mathit{\Omega}_{\text{s}}=\mathcal{O}\setminus \mathit{\Omega}_{\text{f}}$.
Within the design domain $D$, the fluid domain $\mathit{\Omega}_{\text{f}}$ and solid domain $\mathit{\Omega}_{\text{s}}$ are indicated by $\gamma (\mathbf{x})=1$ and $\gamma(\mathbf{x})=0$, respectively, where $\gamma(\mathbf{x})$ is the design variable field of a material distribution.
The evaluation functions are given by the governing equations for fluid velocity $\mathbf{u}: \mathit{\Omega}_{\text{f}}\to\mathbb{R}^2$, pressure $p: \mathit{\Omega}_{\text{f}}\to\mathbb{R}$, and temperature $T: \mathit{\Omega}_{\text{f}}\to\mathbb{R}$.

In this paper, we model the turbulent flow using the steady-state incompressible Reynolds-averaged Navier-Stokes (RANS) equations with the standard $k$-$\epsilon$ model, which is a widely used approximate model for turbulent flow analysis.
The governing equations for the fluid velocity $\mathbf{u}$ and pressure $p$ in $\mathit{\Omega}_{\text{f}}$ are given by the following dimensionless RANS equations:
\begin{equation}\label{eq_cont}
  \nabla\cdot\mathbf{u} = 0,
\end{equation}
\begin{equation}\label{eq_rans}
  \left(\mathbf{u}\cdot\nabla\right)\mathbf{u} = -\nabla p + \nabla\cdot \left(\frac{2}{Re}\mathbf{S}\right)+\nabla\cdot\mathbf{T},
\end{equation}
where $Re$ is the Reynolds number.
$\mathbf{S}=\left(\nabla\mathbf{u}+(\nabla\mathbf{u})^\top\right)/2$ is the mean strain-rate tensor and $\mathbf{T}=2\nu _{\text{t}}\mathbf{S}-2k\mathbf{I}/3$ is the Reynolds stress tensor based on the Boussinesq approximation.
Here, $k$ is the turbulent kinetic energy and $\nu_{\text{t}}$ is the eddy viscosity, which is defined as follows:
\begin{equation}\label{eq_nut}
  \nu_{\text{t}} = C_\mu\frac{k^2}{\epsilon},
\end{equation}
where $C_\mu$ is a empirical constant, and $\epsilon$ is the energy dissipation rate.
$k$ and $\epsilon$ are governed by the following equations:
\begin{equation}\label{eq_k}
  \mathbf{u}\cdot\nabla k = \nabla\cdot\left(\frac{\nu_{\text{t}}}{\sigma_k}\nabla k\right) + 2\nu_{\text{t}}\mathbf{S}:\mathbf{S} - \epsilon,
\end{equation}
\begin{equation}\label{eq_epsilon}
  \mathbf{u}\cdot\nabla\epsilon = \nabla\cdot\left(\frac{\nu_{\text{t}}}{\sigma_\epsilon}\nabla \epsilon\right) + C_{1\epsilon}\frac{\epsilon}{k}2\nu_{\text{t}}\mathbf{S}:\mathbf{S} - C_{2\epsilon}\frac{\epsilon^2}{k},
\end{equation}
where the empirical constants are set as $C_\mu=0.09$, $C_{1\epsilon}=1.44$, $C_{2\epsilon}=1.92$, $\sigma_k=1.0$, and $\sigma_\epsilon=1.3$.
The temperature $T: \mathit{\Omega}_{\text{f}}\to\mathbb{R}$ is governed by the following conjugate heat transfer equation in a dimensionless form:
\begin{equation}\label{eq_heathf}
  \mathbf{u}\cdot\nabla T = \nabla\cdot\left(\left(\frac{1}{Pe}+\frac{\nu_{\text{t}}}{Pr_{\text{t}}}\right)\nabla T\right),
\end{equation}
where $Pe$ is the P\'{e}clet number defined as $Re\cdot Pr$ where $Pr$ is the Prandtl number, and $Pr_{\text{t}}$ is the turbulent Prandtl number.
The parameters regarding these governing equations for the HF model are listed in Table~\ref{tab_thermofluidhf}.

\begin{table}[t]
  \centering
  \caption{Parameters regarding the high-fidelity model in the two-dimensional turbulent heat transfer problem.}\label{tab_thermofluidhf}
  \begin{tabular*}{\textwidth}{@{\extracolsep\fill}lll}
    \toprule
    Description & Symbol & Value \\
    \midrule
    Reynolds number & $Re$ & $5000$ \\
    P\'{e}clet number & $Pe$ & $500$ \\
    Turbulent Prandtl number & $Pr_{\text{t}}$ & $0.1$ \\
    Biot number & $Bi$ & $0.1$ \\
    Filter radius for HF evaluation & $R_{\text{h}}$ & $0.021$ \\
    \bottomrule
  \end{tabular*}
\end{table}

The boundary conditions in the analysis domain shown in Fig.~\ref{problem_thermofluid}(a) are defined in dimensionless forms as follows:
\begin{align}
  \mathbf{u} &= -\mathbf{n}, & T &= 0 && \text{on } \mathit{\Gamma}_{\text{in}}, \\
  p &= 0, & \mathbf{n} \cdot \nabla T &= 0 && \text{on } \mathit{\Gamma}_{\text{out}}, \\
  \mathbf{u}\cdot\mathbf{n} &= 0, & \mathbf{n} \cdot \nabla T &= 0 && \text{on } \mathit{\Gamma}_{\text{wall}},\,\mathit{\Gamma}_{\text{sym}}, \\
  & & -\mathbf{n}\cdot \mathbf{q} &= Bi(1-T) && \text{on } \partial\mathit{\Omega}_{\text{f}}\setminus \left(\mathit{\Gamma}_{\text{wall}}\cup \mathit{\Gamma}_{\text{sym}}\right),
\end{align}
where $\mathbf{n}$ is the outward unit normal vector and $Bi$ is the Biot number.
The heat flux $\mathbf{q}$ in a dimensionless form is defined as follows:
\begin{equation}\label{eq_heatflux}
  \mathbf{q} = -\left(\frac{1}{Pe}+\frac{\nu_{\text{t}}}{Pr_{\text{t}}}\right)\nabla T.
\end{equation}

Based on these governing equations and boundary conditions, the optimization problem is formulated as follows:
\begin{eqnarray}\label{eq_thermofluidhf}
  \begin{aligned}
    & \underset{\gamma(\mathbf{x})}{\text{minimize}}
    && J_1 = -\frac{\int_{\mathit{\Gamma}_{\text{out}}}T\,\mathrm{d}\mathit{\Gamma}}{\int_{\mathit{\Gamma}_{\text{out}}}\mathrm{d}\mathit{\Gamma}}, \\
    &&& J_2 = \frac{\int_{\mathit{\Gamma}_{\text{in}}}P\,\mathrm{d}\mathit{\Gamma}}{\int_{\mathit{\Gamma}_{\text{in}}}\mathrm{d}\mathit{\Gamma}}, \\
    & \text{subject to}
    && \gamma(\mathbf{x})\in\{0, 1\}, \quad \forall\mathbf{x}\in D.
  \end{aligned}
\end{eqnarray}
Herein, setting $J_1$ and $J_2$ as the objective functions corresponds to maximizing the heat exchange amount between the fluid and solid domains, and minimizing the pressure loss in the fluid domain, respectively.

As shown in Fig.~\ref{problem_thermofluid}(c), the analysis domain $\mathcal{O}$ is discretized using body-fitted meshes with triangular elements and five boundary layers with quadrilateral elements, which are automatically built with the maximum and minimum element sizes of $6.6\times 10^{-2}$ and $2.9\times 10^{-3}$, respectively.
All state variable fields are discretized using $\mathbb{P}^1$ Lagrange finite elements.

\begin{table}[t]
  \centering
  \caption{Parameters regarding the low-fidelity model in the two-dimensional turbulent heat transfer problem.}\label{tab_thermofluidlf}
  \begin{tabular*}{\textwidth}{@{\extracolsep\fill}lll}
    \toprule
    Description & Symbol & Value \\
    \midrule
    Inverse permeability & $\alpha$ & $1.0\times 10^4$ \\
    Prandtl number of fluid & $Pr_{\text{f}}$ & $0.1$ \\
    Prandtl number of solid & $Pr_{\text{s}}$ & $1$ \\
    Tuning parameter for inverse permeability & $q_{\alpha}$ & $100$ \\
    Tuning parameter for Prandtl number & $q_{Pr}$ & $10$ \\
    Minimum Reynolds number & $Re_{\text{min}}$ & $10$ \\
    Maximum Reynolds number & $Re_{\text{max}}$ & $80$ \\
    Minimum volumetric heat transfer coefficient & $\beta_{\text{min}}$ & $0.1$ \\
    Maximum volumetric heat transfer coefficient & $\beta_{\text{max}}$ & $0.5$ \\
    Number of seeding parameters for Reynolds number & $n_{s_1}$ & 4 \\
    Number of seeding parameters for volumetric heat transfer coefficient & $n_{s_2}$ & 25 \\
    Filter radius & $R$ & $0.06$ \\
    \bottomrule
  \end{tabular*}
\end{table}

\subsubsection*{Low-fidelity model}
The original problem based on the HF model is formulated with a turbulent flow model under high Reynolds number conditions, which is intractable for gradient-based topology optimization due to its strong nonlinearity.
In order to stably solve the LF optimization problem and derive diverse and promising initial designs, we employ a laminar flow model under low Reynolds number conditions as the LF model, which is a common approach in the literature~\cite{alexandersen2020}.
To enable density-based topology optimization by representing both fluid and solid domains using a single design variable, we consider the following governing equations for the state variable fields $\tilde{\mathbf{u}}: \mathcal{O}\to\mathbb{R}^2$, $\tilde{p}: \mathcal{O}\to\mathbb{R}$, and $\tilde{T}: \mathcal{O}\to\mathbb{R}$ in the LF model:
\begin{equation}\label{eq_contlf}
  \nabla\cdot\tilde{\mathbf{u}} = 0,
\end{equation}
\begin{equation}\label{eq_nslf}
  \left(\tilde{\mathbf{u}}\cdot\nabla\right)\tilde{\mathbf{u}} = -\nabla \tilde{p} + \frac{1}{\tilde{Re}}\nabla^2\tilde{\mathbf{u}} - \alpha_{\tilde{\gamma}} \tilde{\mathbf{u}},
\end{equation}
\begin{equation}\label{eq_heathflf}
  \tilde{\mathbf{u}}\cdot\nabla \tilde{T} = \frac{1}{\tilde{Re}\tilde{Pr}_{\tilde{\gamma}}}\nabla^2\tilde{T} + \beta_{\tilde{\gamma}}(1-\tilde{T}),
\end{equation}
where $\tilde{Re}$ is the Reynolds number for the laminar flow model.
Note that the governing equations in Eqs.~\eqref{eq_contlf}, \eqref{eq_nslf}, and \eqref{eq_heathflf} are formulated in dimensionless forms.
The design variable-dependent parameters $\alpha_{\tilde{\gamma}}$, $\beta_{\tilde{\gamma}}$, and $\tilde{Pr}_{\tilde{\gamma}}$ are defined as follows:
\begin{equation}\label{eq_alpha}
  \alpha_{\tilde{\gamma}} = \alpha\frac{1-\tilde{\gamma}}{\tilde{\gamma}+q_{\alpha}},
\end{equation}
\begin{equation}\label{eq_beta}
  \beta_{\tilde{\gamma}} = \beta(1-\tilde{\gamma}),
\end{equation}
\begin{equation}\label{eq_pr}
  \tilde{Pr}_{\tilde{\gamma}} = Pr_{\text{f}}+(Pr_{\text{s}}-Pr_{\text{f}})\frac{1-\tilde{\gamma}}{1+q_{Pr}\tilde{\gamma}},
\end{equation}
where $\tilde{\gamma} = \sum_{e=1}^n \gamma_e^{(k)} \chi_e(\mathbf{x})$ is the design variable field converted from $\boldsymbol{\gamma}^{(k)}$ to a continuous field, where $\chi_e(\mathbf{x})$ is the indicator function of the $e$-th element, i.e., $\chi_e(\mathbf{x})=1$ if $\mathbf{x}$ is in the $e$-th element and $\chi_e(\mathbf{x})=0$ otherwise.
The constants $\alpha$ and $\beta$ denote the inverse permeability and volumetric heat transfer coefficient, respectively.
$Pr_{\text{f}}$ and $Pr_{\text{s}}$ are the Prandtl numbers of the fluid and solid domains, respectively, and $q_{\alpha}$ and $q_{Pr}$ are tuning parameters controlling the convexity of the interpolation functions.
$\alpha_{\tilde{\gamma}}$ appears in the Navier-Stokes equation in Eq.~\eqref{eq_nslf} to model the solid domain~\cite{borrvall2003}, while $\beta_{\tilde{\gamma}}$ is used in the energy equation in Eq.~\eqref{eq_heathflf} to represent the heat source in the solid domain.

\begin{figure}[t]
  \centering
  \includegraphics[width=.8\textwidth]{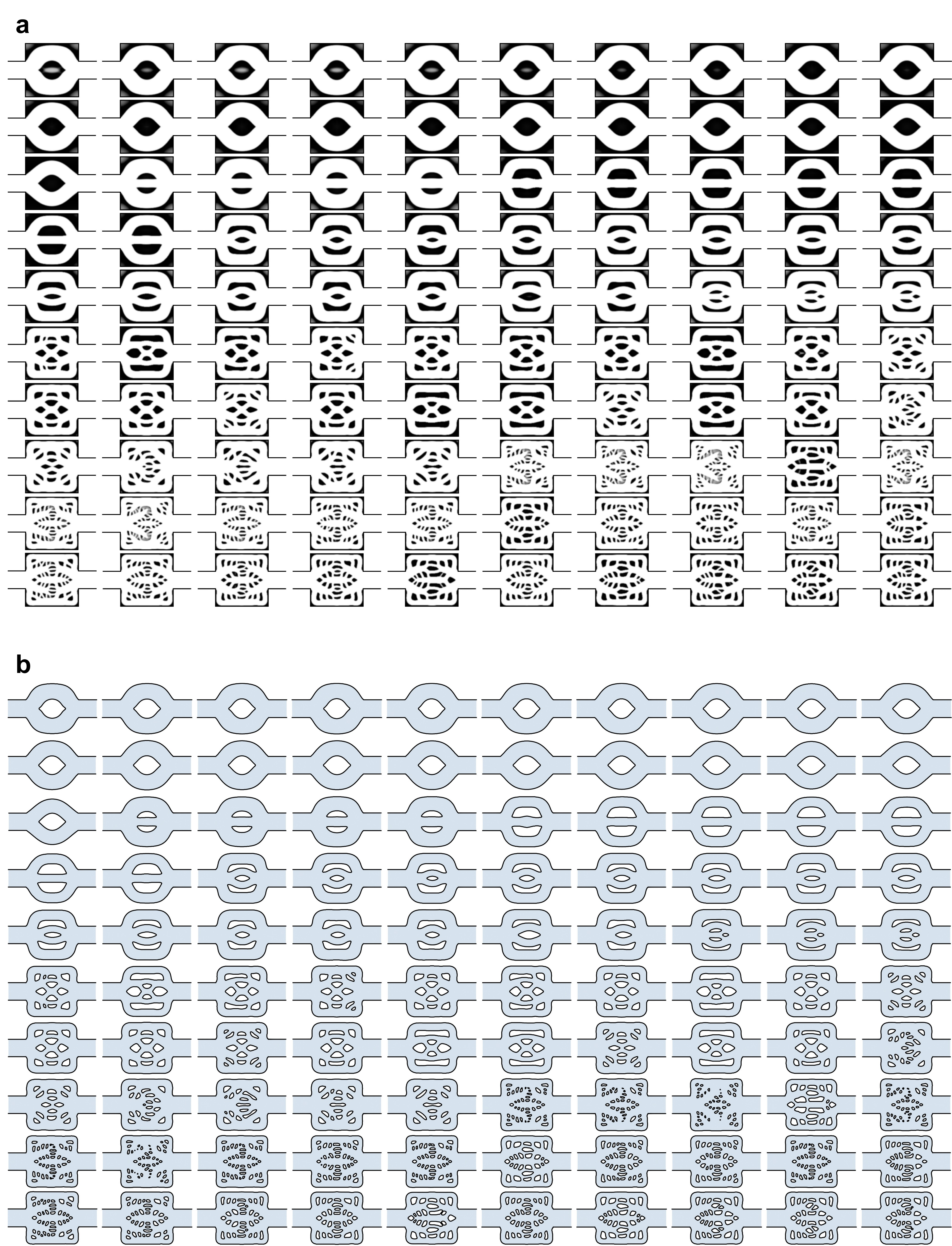}
  \caption{Low-fidelity optimized initial designs in the two-dimensional heat sink design problem: (a) density distributions on the low-fidelity model; (b) material distributions on the high-fidelity model.}\label{result_thermofluid_ini}
\end{figure}

The boundary conditions in the analysis domain shown in Fig.~\ref{problem_thermofluid}(a) are defined as follows:
\begin{align}
  \tilde{p} &= 1, & \tilde{T} &= 0 && \text{on } \mathit{\Gamma}_{\text{in}}, \\
  \tilde{p} &= 0, & \mathbf{n}\cdot\nabla \tilde{T} &= 0 && \text{on } \mathit{\Gamma}_{\text{out}}, \\
  \tilde{\mathbf{u}} &= 0, & \mathbf{n}\cdot\nabla \tilde{T} &= 0 && \text{on } \mathit{\Gamma}_{\text{wall}}, \\
  \tilde{\mathbf{u}}\cdot\mathbf{n} &= 0, & \mathbf{n}\cdot\nabla \tilde{T} &= 0 && \text{on } \mathit{\Gamma}_{\text{sym}}. 
\end{align}

\begin{figure}[t]
  \centering
  \includegraphics[width=.95\textwidth]{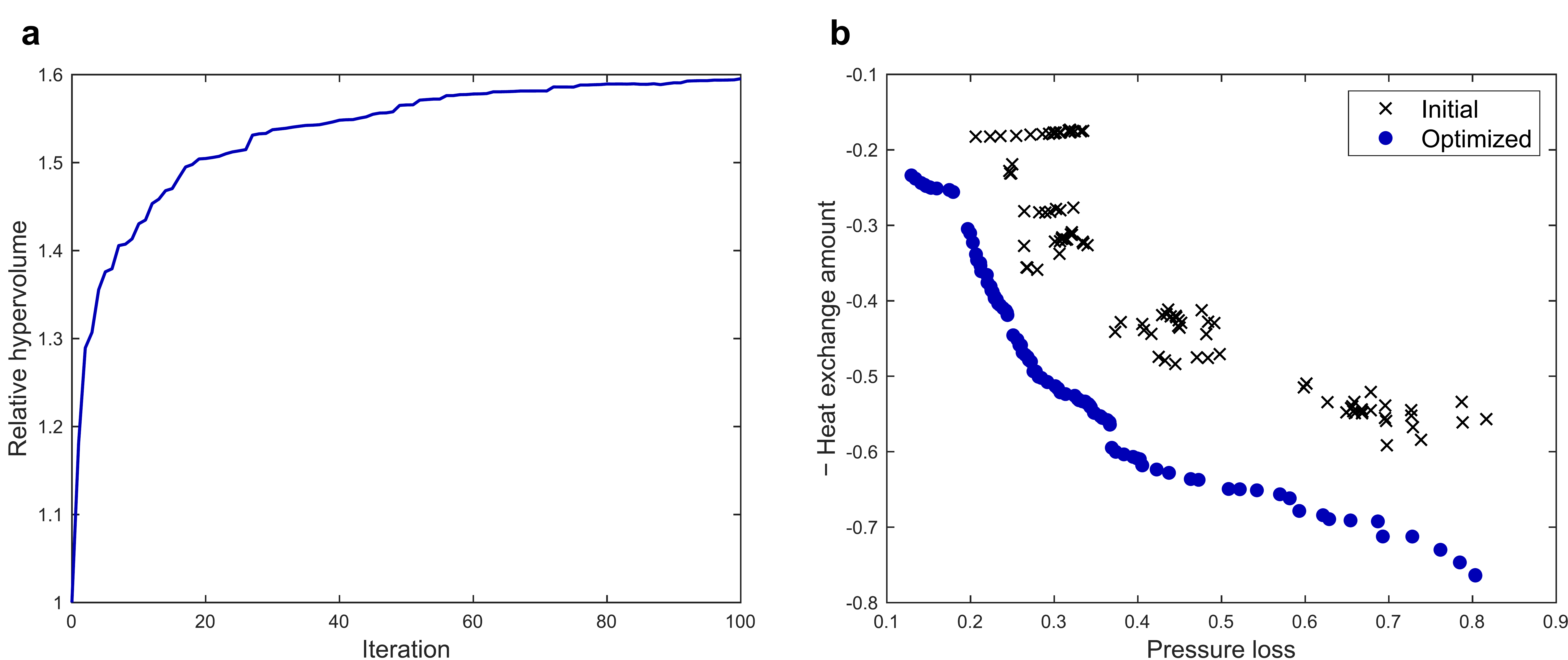}
  \caption{Optimization results in the two-dimensional heat sink design problem: (a) convergence history of the hypervolume indicator normalized by the initial value; (b) objective space where the performance values of the initial and optimized designs are plotted.}\label{result_thermofluid}
\end{figure}

\begin{figure}[t]
  \centering
  \includegraphics[width=.8\textwidth]{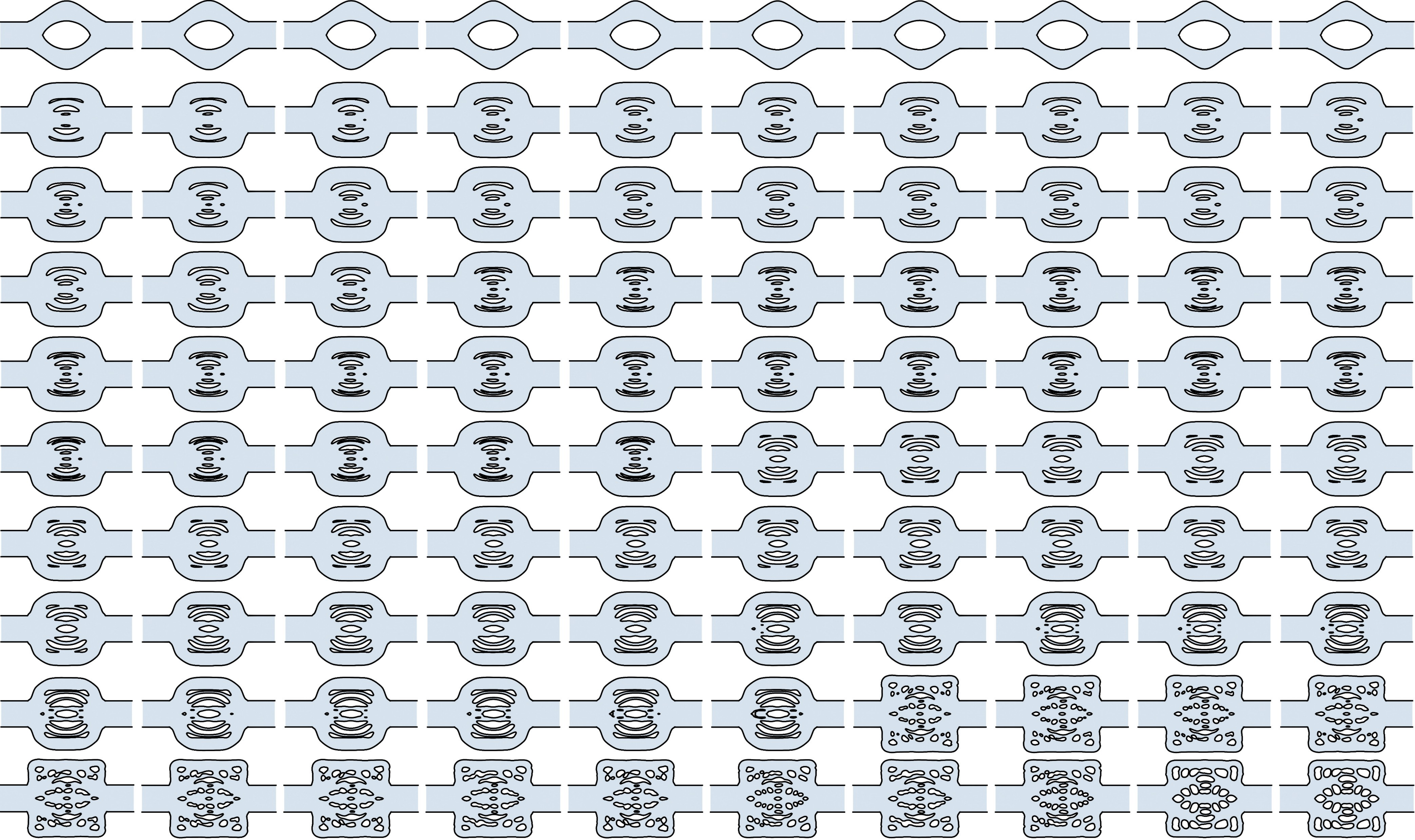}
  \caption{Optimized designs in the two-dimensional heat sink design problem}\label{result_thermofluid_opt}
\end{figure}

Based on these governing equations and boundary conditions, the LF optimization problem is formulated as follows:
\begin{eqnarray}\label{eq_thermofluidlf}
  \begin{aligned}
    & \underset{\boldsymbol{\gamma}^{(k)}}{\text{minimize}}
    && \widetilde{J} = -\int_{\mathit{\Gamma}_{\text{out}}}(\tilde{\mathbf{u}}\cdot\mathbf{n})\tilde{T}\,\mathrm{d}\mathit{\Gamma} \\
    & \text{subject to}
    && \boldsymbol{\gamma}^{(k)}\in [0, 1]^n, \\
    & \text{for given } && \mathbf{s}^{(k)}=[s_1, s_2]^\top.
  \end{aligned}
\end{eqnarray}
Herein, two seeding parameters, $s_1$ and $s_2$, are employed to seed the Reynolds number $\tilde{Re}$ in Eqs.~\eqref{eq_nslf} and \eqref{eq_heathflf}, and the volumetric heat transfer coefficient $\beta$ in Eq.~\eqref{eq_beta}, respectively.
The Reynolds number $\tilde{Re}$ is obtained by linearly interpolating between $Re_{\text{min}}$ and $Re_{\text{max}}$ according to $s_1$, while $\beta$ is interpolated between $\beta_{\text{min}}$ and $\beta_{\text{max}}$ using $s_2$.
Both $s_1$ and $s_2$ are sampled uniformly from the interval $[0,1]$ with $n_{s_1}$ and $n_{s_2}$ divisions, respectively, as in the numerical example in Section~\ref{sec51}.
Consequently, the LF optimization problem of Eq.~\eqref{eq_thermofluidlf} is solved independently $N_{\text{lf}} = n_{s_1} \times n_{s_2}$ times.
As shown in Fig.~\ref{problem_thermofluid}(b), the design variable $\boldsymbol{\gamma}^{(k)}$ is discretized using structured meshes with 10,368 quadrilateral elements, i.e., $n=10368$.
The overall parameters for the LF model are listed in Table~\ref{tab_thermofluidlf}.

\begin{figure}[t]
  \centering
  \includegraphics[width=\textwidth]{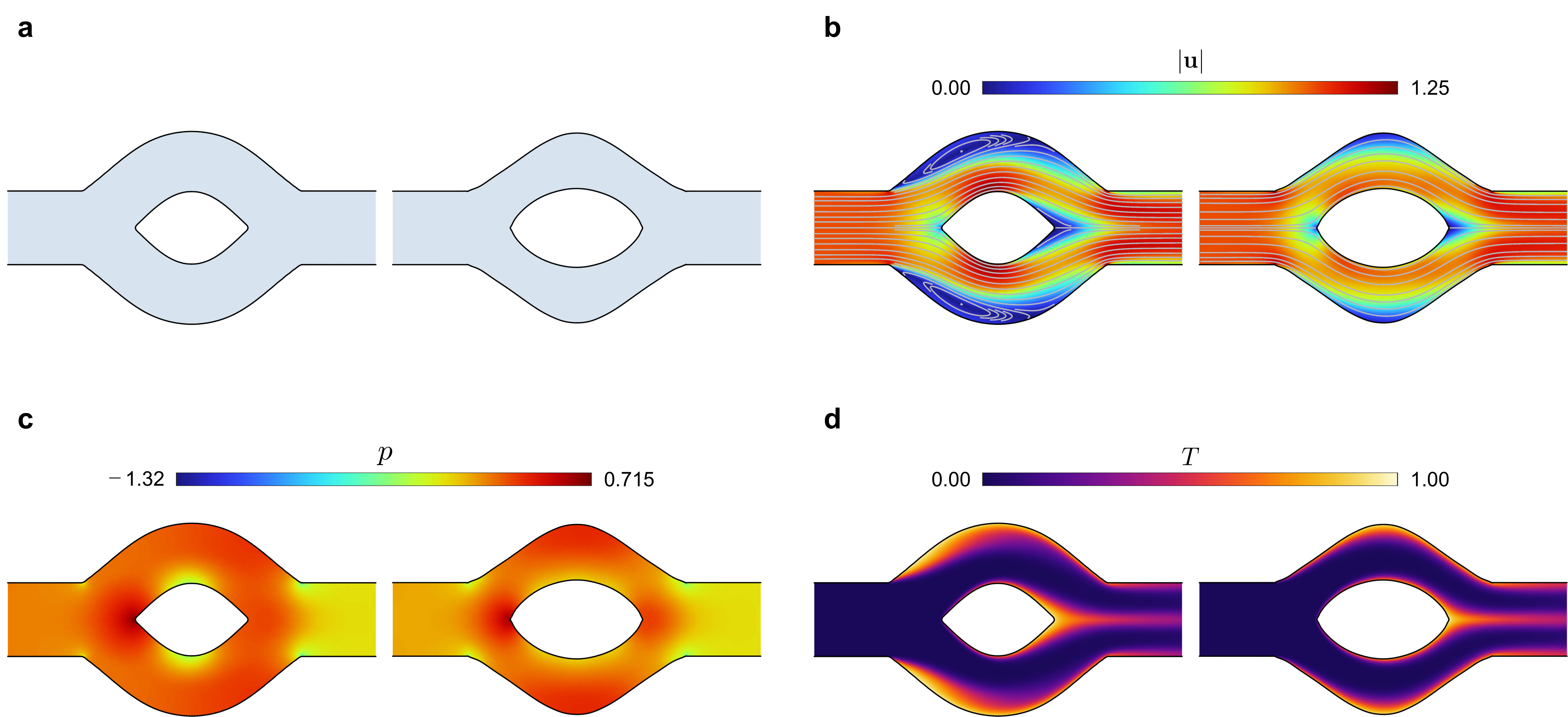}
  \caption{Results comparison of initial (left) and optimized (right) designs with similar shapes in the two-dimensional heat sink design problem: (a) channel configurations; (b) velocity magnitude $|\mathbf{u}|$ with stream lines; (c) pressure $p$; (d) temperature $T$. Herein, these objective values are: $J_1=-0.183$, $J_2=0.206$ (left); $J_1=-0.234$, $J_2=0.129$ (right), where $J_1$ and $J_2$ are the negative heat exchange amount and the pressure loss, respectively.}\label{result_thermofluid_compare1}
\end{figure}

\begin{figure}[t]
  \centering
  \includegraphics[width=\textwidth]{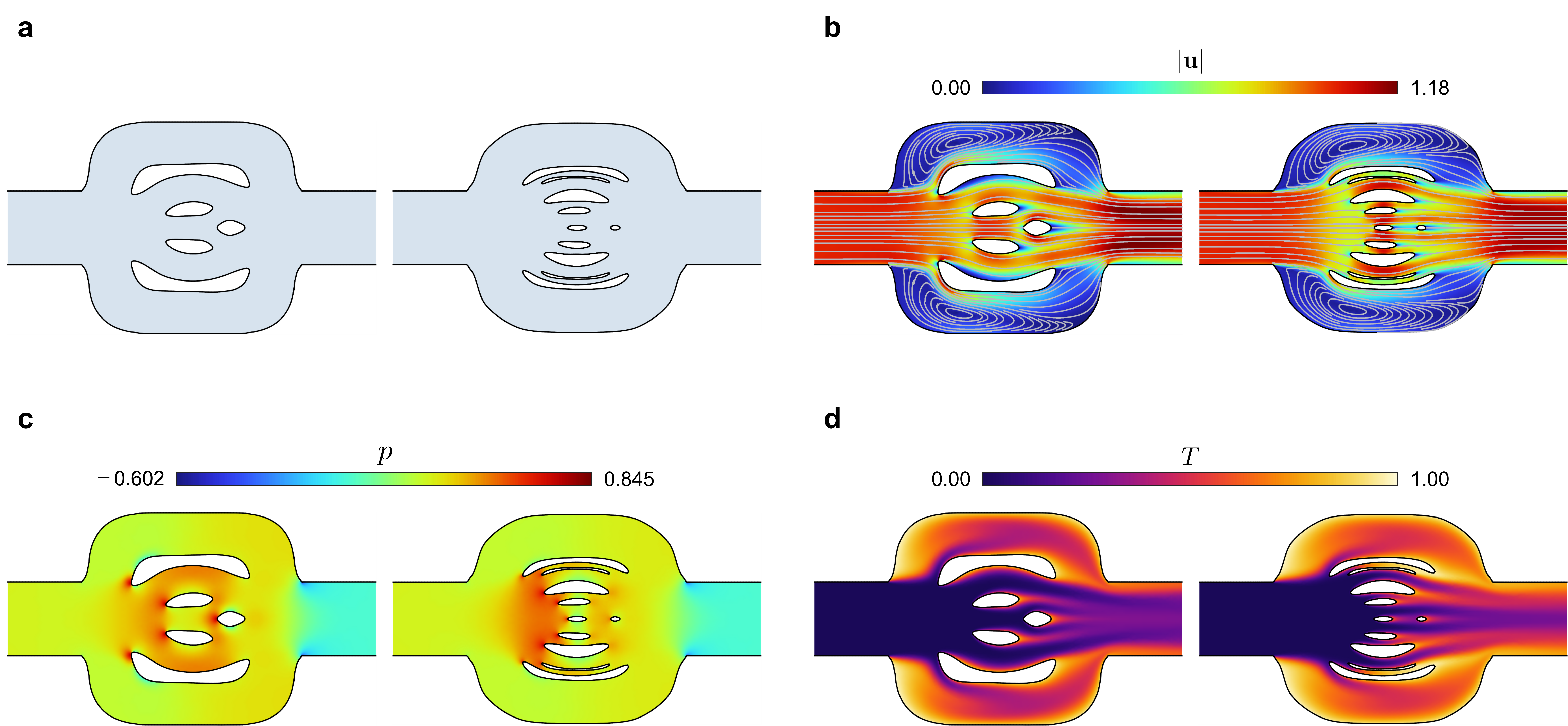}
  \caption{Results comparison of initial (left) and optimized (right) designs under nearly identical pressure loss conditions in the two-dimensional heat sink design problem: (a) channel configurations; (b) velocity magnitude $|\mathbf{u}|$ with stream lines; (c) pressure $p$; (d) temperature $T$. Herein, these objective values are: $J_1=-0.359$, $J_2=0.280$ (left); $J_1=-0.501$, $J_2=0.281$ (right), where $J_1$ and $J_2$ are the negative heat exchange amount and the pressure loss, respectively.}\label{result_thermofluid_compare2}
\end{figure}

\begin{figure}[t]
  \centering
  \includegraphics[width=\textwidth]{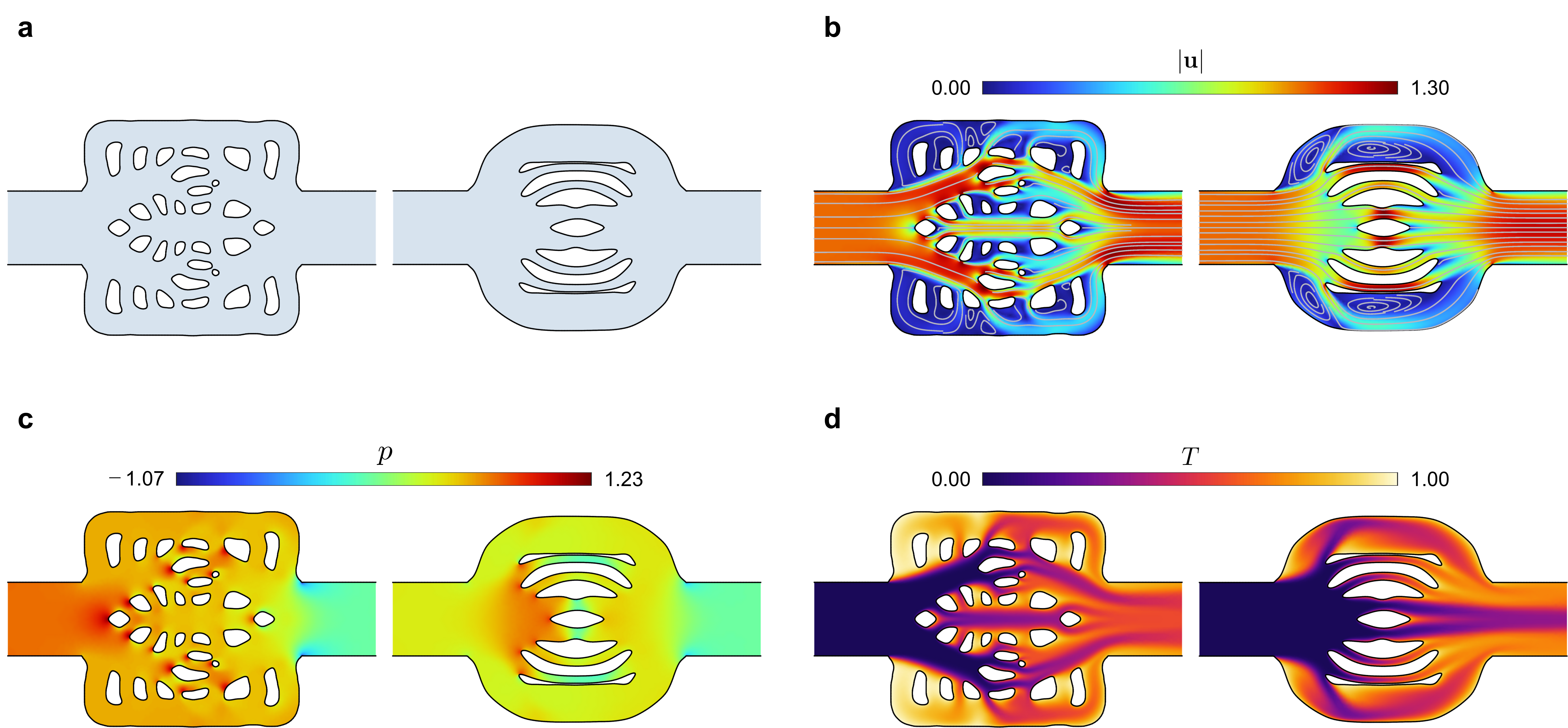}
  \caption{Results comparison of initial (left) and optimized (right) designs under nearly identical heat exchange amount conditions in the two-dimensional heat sink design problem: (a) channel configurations; (b) velocity magnitude $|\mathbf{u}|$ with stream lines; (c) pressure $p$; (d) temperature $T$. Herein, these objective values are: $J_1=-0.591$, $J_2=0.670$ (left); $J_1=-0.595$, $J_2=0.269$ (right), where $J_1$ and $J_2$ are the negative heat exchange amount and the pressure loss, respectively.}\label{result_thermofluid_compare3}
\end{figure}

\subsubsection{Results and discussion}\label{sec522}

Figure~\ref{result_thermofluid_ini} shows the initial population obtained by solving the LF optimization problem of Eq.~\eqref{eq_thermofluidlf}.
As in the example in Section~\ref{sec51}, seeding parameters $s_1$ and $s_2$ in Eqs.~\eqref{eq_thermofluidlf} yield diverse density distributions on the LF model shown in Fig.~\ref{result_thermofluid_ini}(a), which are then appropriately converted into flow channels on the HF model as shown in Fig.~\ref{result_thermofluid_ini}(b).

Figure~\ref{result_thermofluid} shows the optimization results, including the hypervolume convergence history and the objective space.
In Fig.~\ref{result_thermofluid}(a), the steady increase in hypervolume over iterations indicates that the initial designs are progressively evolved by Wasserstein crossover, finally achieving a 60\% improvement from the initial value.
As shown in Fig.~\ref{result_thermofluid}(b), both objective values have significantly improved from the initial designs, and the optimized solutions form a well-aligned Pareto front.

Figure~\ref{result_thermofluid_opt} shows the optimized population.
Compared with the initial population in Fig.~\ref{result_thermofluid_ini}(b), the optimized population includes both designs that resemble the initial ones and others with entirely novel channel configurations.
To further investigate their differences, Fig.~\ref{result_thermofluid_compare1} compares the initial and optimized channel configurations that share similar overall shapes.
A particularly notable difference is in the streamlines shown in Fig.~\ref{result_thermofluid_compare1}(b): while the initial design exhibits characteristic turbulent vortices near the top and bottom ends of the channel, the optimized design suppresses such vortices by narrowing these regions.
As a result, there is a significant difference in pressure loss $J_2$ between them.

Figure~\ref{result_thermofluid_compare2} compares the initial and optimized designs under nearly identical pressure loss conditions.
Focusing on the solid domains given as heat sources in Fig.~\ref{result_thermofluid_compare2}(a), the initial design contains five such domains, whereas the optimized one has ten.
As a result, the optimized configuration achieves approximately 40\% better heat exchange amount $J_1$ than the initial one.
Nevertheless, as shown in Fig.~\ref{result_thermofluid_compare2}(b), the pressure loss $J_2$ remains comparable between them by arranging solid domains so as not to obstruct the flow.

Figure~\ref{result_thermofluid_compare3} presents a comparison under approximately identical heat exchange amount conditions.
As shown in Fig.~\ref{result_thermofluid_compare3}(a), the initial and optimized designs exhibit entirely different channel configurations: the initial one resembles a pin-fin heat sink, while the optimized one takes on a plate-fin-like configuration with an adjusted outer shape.
This unique curved plate-fin configuration reduces the pressure loss $J_2$ by approximately 60\% compared to the initial design.
In Fig.~\ref{result_thermofluid_compare3}(b), the initial flow channel, originally optimized under a laminar flow model, shows minimal flow near the solid domains close to the inlet under the turbulent flow model.
In contrast, the optimized flow channel maintains high flow velocity between all fins, indicating that effective heat exchange has been achieved with fewer solid domains than the initial one.

These results indicate the potential of the proposed framework to handle complex physical models such as turbulent flow.
In particular, the emergence of novel flow channel configurations---absent from the initial population but introduced through Wasserstein crossover during the optimization process---appears to be a key factor in achieving significant performance improvements.

\begin{figure}[t]
  \centering
  \includegraphics[width=.75\textwidth]{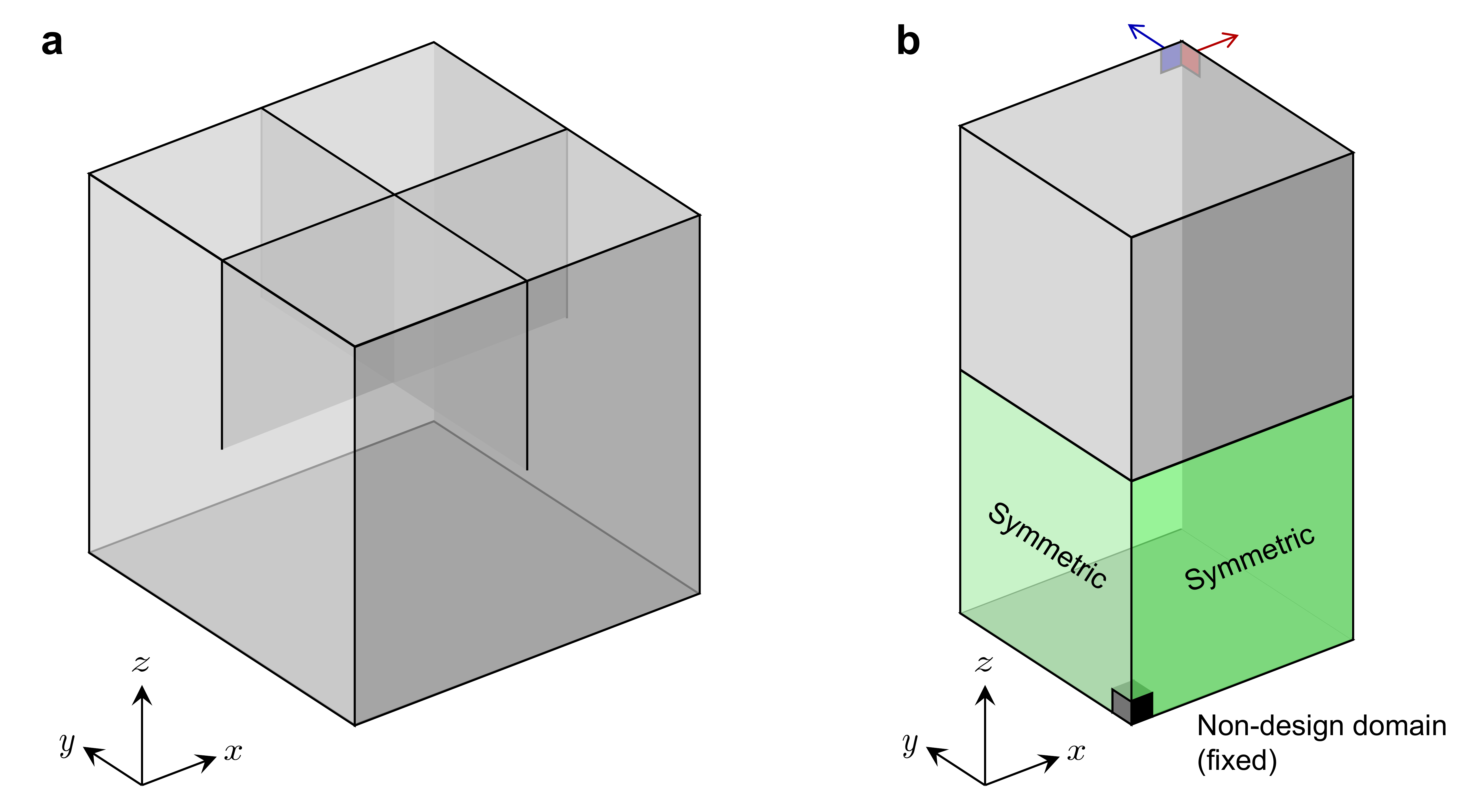}
  \caption{Problem settings in the three-dimensional cracked box design problem: (a) the whole domain; (b) the design domain and boundary conditions defined in the quarter domain.}\label{problem_3dcrack}
\end{figure}

\begin{table}[t]
  \centering
  \caption{Overall parameters in the three-dimensional stress minimization problem.}\label{tab_3dcrack}
  \begin{tabular*}{\textwidth}{@{\extracolsep\fill}lll}
    \toprule
    Description & Symbol & Value \\
    \midrule
    Maximum number of iterations & $t_{\text{max}}$ & 150 \\
    Population size & $N_{\text{pop}}$ & 60 \\ 
    Number of offspring & $N_{\text{xo}}$ & 60 \\
    Number of LF optimized initial designs & $N_{\text{lf}}$ & 120 \\
    Number of seeding parameters for filter radius & $n_{s_1}$ & 3 \\
    Number of seeding parameters for volume fraction & $n_{s_2}$ & 20 \\
    Number of seeding parameters for norm parameter & $n_{s_3}$ & 2 \\
    Minimum filter radius & $R_{\text{min}}$ & $0.04$ \\
    Maximum filter radius & $R_{\text{max}}$ & $0.12$ \\
    Minimum volume fraction & $V_{\text{min}}$ & $0.10$ \\
    Maximum volume fraction & $V_{\text{max}}$ & $0.30$ \\
    Minimum norm parameter for LF model & $P_{\text{min}}$ & 6 \\
    Maximum norm parameter for LF model & $P_{\text{max}}$ & 8 \\
    Filter radius for HF evaluation & $R_{\text{h}}$ & $0.025$ \\
    Minimum entropic regularization coefficient & $\varepsilon_{\min}$ & $1.0\times 10^{-4}$ \\
    Maximum entropic regularization coefficient & $\varepsilon_{\max}$ & $1.0\times 10^{-3}$ \\
    Convergence tolerance in Sinkhorn algorithm & $\tau$ & $1.0\times 10^{-5}$ \\
    \bottomrule
  \end{tabular*}
\end{table}

\begin{figure}[t]
  \centering
  \includegraphics[width=.95\textwidth]{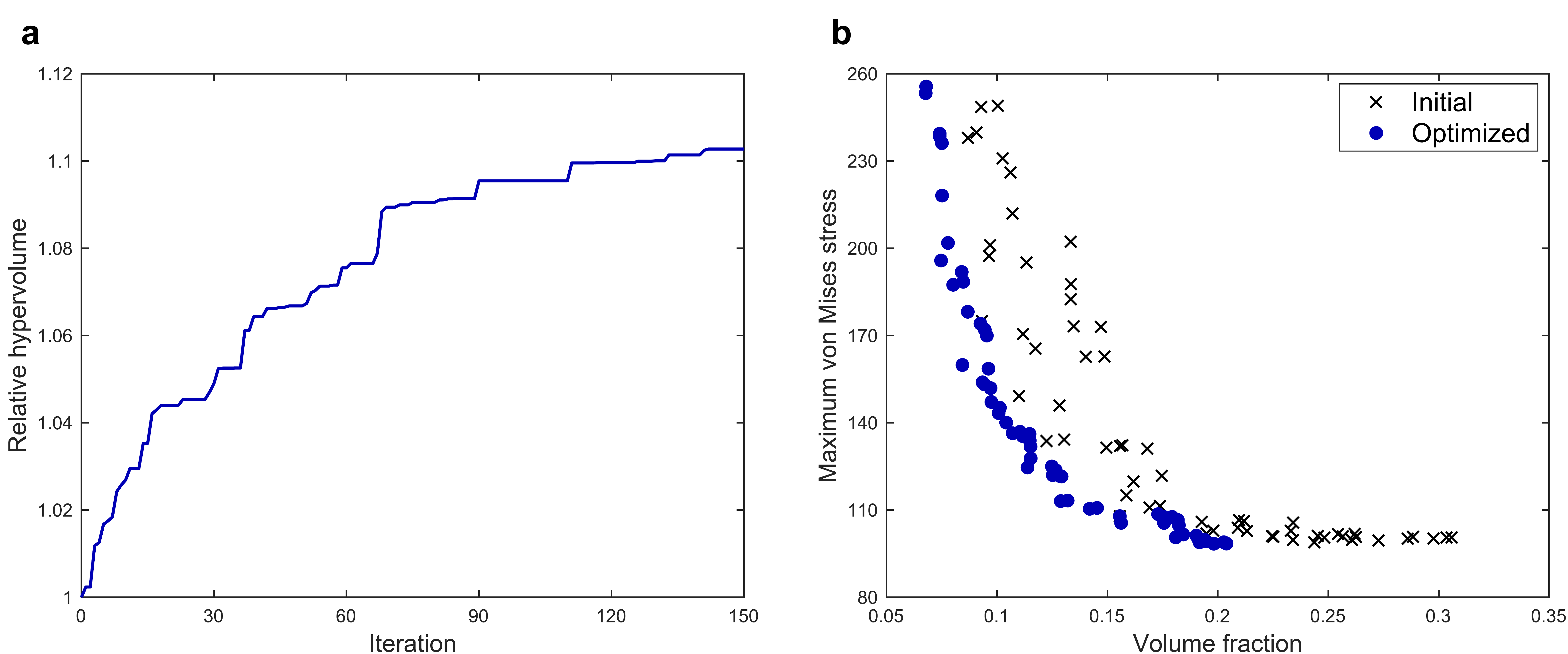}
  \caption{Optimization results in the three-dimensional cracked box design problem: (a) convergence history of the hypervolume indicator normalized by the initial value; (b) objective space where the performance values of the initial and optimized designs are plotted.}\label{result_3dcrack}
\end{figure}

\subsection{3D stress minimization problem}\label{sec53}

\subsubsection{Problem settings}\label{sec531}

Let us consider a three-dimensional stress minimization problem.
Figure~\ref{problem_3dcrack} illustrates the analysis domain of the three-dimensional cracked box.
In Fig.~\ref{problem_3dcrack}(a), two cracks are present in both the $x$ and $y$ directions, and the design domain $D\subset\mathbb{R}^3$ is defined as a quarter of the whole domain shown in Fig.~\ref{problem_3dcrack}(b), where the cracks are represented by applying symmetric boundary conditions on the bottom part.
This problem can be interpreted as a three-dimensional extension of the cracked plate design problem in Section~\ref{sec51}.
The formulations of HF and LF models are basically identical to Eqs.~\eqref{eq_crackhf} and \eqref{eq_cracklf}, respectively, except that an additional seeding parameter $s_3$ is introduced to change the norm parameter $P$ in the $P$-norm stress in Eq.~\eqref{eq_cracklf}.
The norm parameter $P$ is obtained by linearly interpolating between $P_{\text{min}}$ and $P_{\text{max}}$ according to $s_3$, where $P_{\text{min}}$ and $P_{\text{max}}$ are the minimum and maximum values of the norm parameter $P$, respectively.
The parameters regarding the overall algorithm are listed in Table~\ref{tab_3dcrack}.

In the LF model, the design variable field $\boldsymbol{\gamma}^{(k)}$ is discretized using structured meshes with $40\times 40\times 80$ hexahedral elements, i.e., $n=128000$, whereas the HF model employs body-fitted meshes with tetrahedral elements whose maximum and minimum element sizes are set to $0.11$ and $8.5\times 10^{-3}$, respectively.
As in the two-dimensional case, the design variable $\boldsymbol{\gamma}^{(k)}$ is discretized using $\mathbb{P}^0$ Lagrange finite elements for the LF model, while the HF model employs $\mathbb{P}^2$ Lagrange finite elements for structural analysis.

\begin{figure}[t]
  \centering
  \includegraphics[width=.975\textwidth]{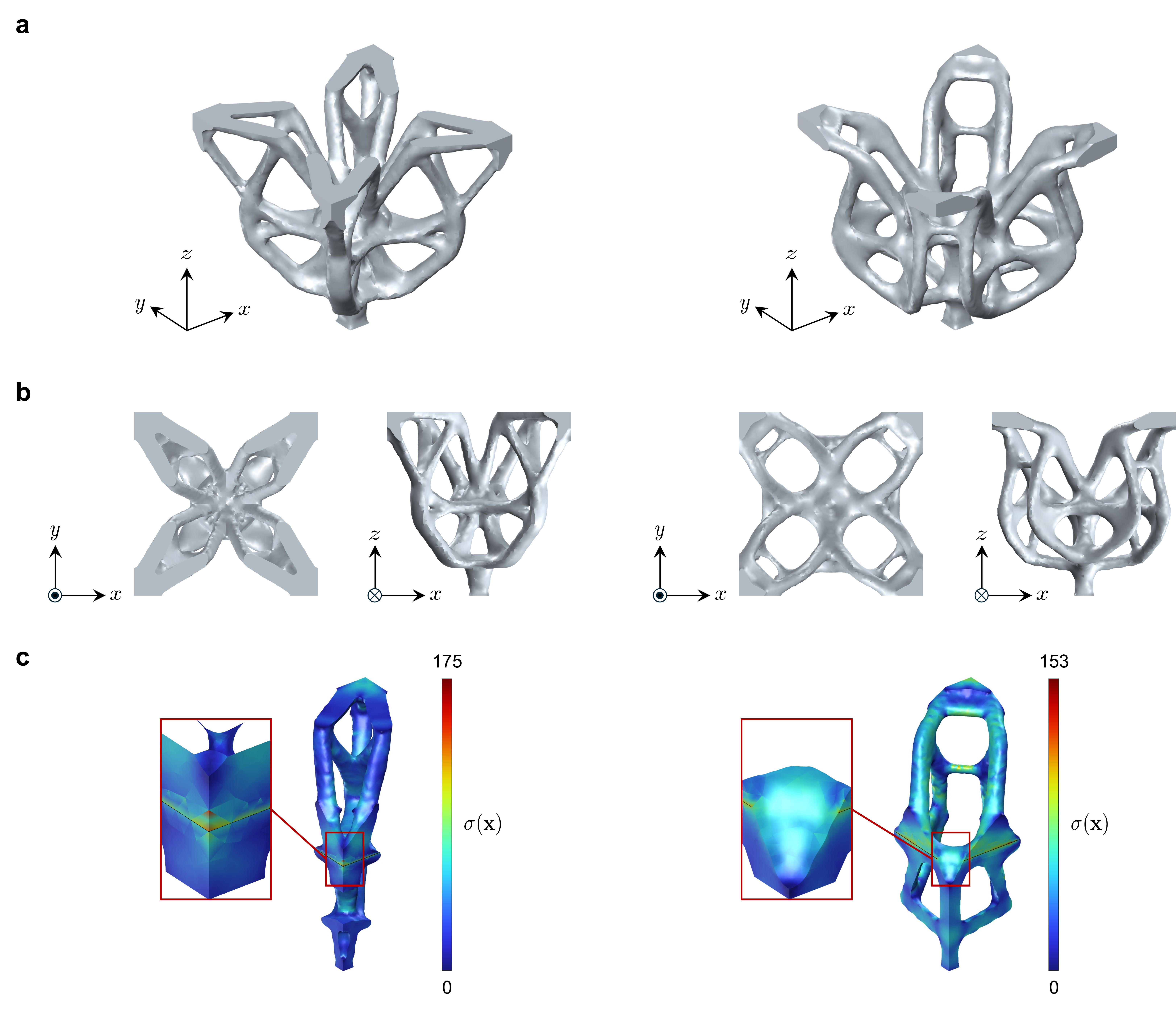}
  \caption{Results comparison of initial (left) and optimized (right) designs under nearly identical volume fraction conditions in the three-dimensional cracked box design problem: (a) material distributions; (b) projections from top and side views; (c) stress distributions $\sigma(\mathbf{x})$. Herein, these objective values are: $J_1=175$, $J_2=9.30\times 10^{-2}$ (left); $J_1=153$, $J_2=9.43\times 10^{-2}$ (right), where $J_1$ and $J_2$ are the maximum stress and volume fraction, respectively.}\label{result_3dcrack_compare}
\end{figure}

\subsubsection{Results and discussion}\label{sec532}

Figure~\ref{result_3dcrack} shows the hypervolume convergence history and the objective space as the optimization results.
The hypervolume improves over iterations in Fig.~\ref{result_3dcrack}(a), finally increasing by about 10\% from its initial value.
In Fig.~\ref{result_3dcrack}(b), the Pareto front has also progressed from the initial solutions.
These results suggest that even in three-dimensional cases, Wasserstein crossover can potentially lead to optimized designs with improved physical performance.

As in the previous numerical examples, we compare the representative initial and optimized designs under nearly identical volume fraction conditions in Fig.~\ref{result_3dcrack_compare}.
The optimized designs shown in Fig.~\ref {result_3dcrack_compare}(a) and (b) are structurally characterized by smoother geometries that include rounded corners and edges, in contrast to the initial designs.
Focusing on the stress distributions in Fig.~\ref{result_3dcrack_compare}(c), it can be seen that the optimized design removes the material at the intersection point of the two cracks, where stress concentration occurs in the initial design.
As a result, the maximum stress $J_1$ is reduced by approximately 13\% in the optimized design compared to the initial one, while the volume fraction $J_2$ remains nearly unchanged.

These results indicate the potential applicability of the proposed framework to three-dimensional topology optimization problems.
As with the two-dimensional cases in Sections~\ref{sec51} and \ref{sec52}, the proposed Wasserstein crossover seems to produce physically reasonable material distributions that may lead to improved performance.

\subsection{Computational cost analysis}\label{sec54}

To assess the scalability of the proposed framework, we measured the computational time for each numerical example.
Table~\ref{tab_time} summarizes the breakdown across the procedures described in Section~\ref{sec3}.
First, for the two-dimensional stress minimization problem in Section~\ref{sec51}, the proposed framework required 31.74 hours, whereas the DDTD framework took 16.84 hours, indicating that Wasserstein crossover incurred approximately 1.9 times higher computational cost compared to VAE-based crossover.
This result suggests that, under the conditions of Section~\ref{sec51}, repeatedly computing the Wasserstein barycenter is more computationally expensive than the training and sampling process of the VAE.
Nevertheless, in light of the striking difference in search performance, as demonstrated by the approximately 2.6 times improvement in hypervolume shown in Fig.~\ref{result_crack}(a), the increased computational cost of Wasserstein crossover appears to be acceptable.
It should also be noted that further acceleration can be expected through parallel computing techniques with multiple GPUs, as the computation of the Wasserstein barycenter for each offspring is independent.

\begin{table}[t]
  \centering
  \caption{Breakdown of computational time for each numerical example (hours).}
  \label{tab_time}
  \begin{tabular*}{\textwidth}{@{\extracolsep\fill}lccc}
    \toprule
    Case & 2D stress minimization & 2D turbulent heat transfer & 3D stress minimization \\
    \midrule
    LF optimization & 0.39 & 0.46 & 19.01 \\
    HF evaluation & 6.36 & 19.58 & 80.02 \\
    Selection & 3.23 & 1.18 & 7.51 \\
    Wasserstein crossover & 21.76 & 21.52 & 24.54 \\
    \midrule
    Total & 31.74 & 42.74 & 131.08 \\
    \bottomrule
  \end{tabular*}
\end{table}

\begin{figure}[t]
  \centering
  \includegraphics[width=.875\textwidth]{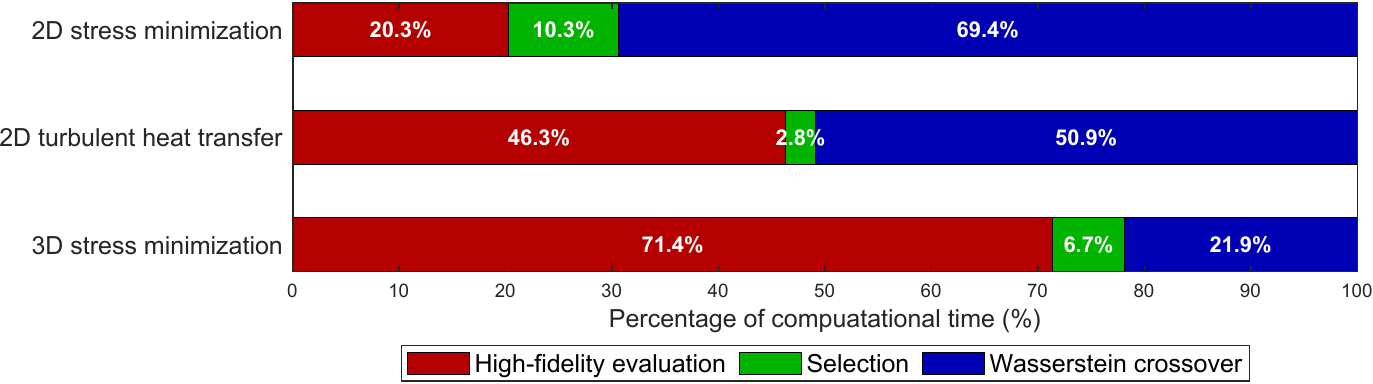}
  \caption{Percentage of computational time within the main loop of the proposed framework for three numerical examples.}\label{result_time}
\end{figure}

Next, comparing the computational times across the three numerical examples in Table~\ref{tab_time} reveals that the total costs are highly dependent on LF optimization and HF evaluation.
In other words, the computational load of the finite element analysis significantly influences the overall cost of the proposed framework.
Regarding the Wasserstein crossover, its computational complexity primarily depends on the number of design variables $n$ and the number of offspring $N_{\text{xo}}$.
Consequently, the computational times for Wasserstein crossover in the two-dimensional problems are nearly identical.
For the three-dimensional stress minimization problem, despite the quite larger number of design variables, the computational time for Wasserstein crossover is not substantially higher, which can be partly attributed to the smaller number of offspring $N_{\text{xo}}$ compared to the two-dimensional cases.
In addition, the stronger entropic regularization achieved by larger values of $\varepsilon_{\min}$ and $\varepsilon_{\max}$, together with a relaxed convergence tolerance $\tau$ in the Sinkhorn algorithm in Algorithm~\ref{alg_barycenter}, led to decrease the number of Sinkhorn iterations, thereby reducing the computational cost.

For a clearer comparison, Fig.~\ref{result_time} illustrates the proportion of computational time spent on each procedure within the main loop, consisting of HF evaluation, selection, and Wasserstein crossover, excluding LF optimization with relatively low computational costs.
In all cases, the majority of the computational time was consumed by HF evaluation and Wasserstein crossover, while the time for selection was negligible.
Notably, even under the identical conditions of parallel computation using half the CPU cores per offspring $N_{\text{xo}}$, its computational cost within the main loop becomes increasingly dominant in problems where the finite element analysis is computationally expensive.
For extending the proposed framework to more complex problems with costly finite element analyses, it will be crucial to reduce the number of evaluations through adjustments in population size or the use of HF surrogate models, which will be explored in future work.

\section{Conclusions}\label{sec6}

In this paper, we introduced a Wasserstein crossover operator for material distributions in topology optimization and developed an EA-based framework.
In the proposed framework, material distributions are represented as density fields following a density-based topology optimization approach.
During the Wasserstein crossover, selected parent candidate solutions are interpreted as probability distributions, and their Wasserstein barycenter is computed to generate offspring.
This approach allows for smooth and meaningful interpolation between parents based on the optimal transport theory, leading to a more reasonable crossover operator compared to conventional ones, such as the VAE-based crossover employed in DDTD.
The proposed framework also incorporates the idea of multifidelity design: LF optimization, which is formulated as an easily solvable pseudo-problem, is used to derive initial designs, followed by iterative cycles of HF evaluation, selection, and Wasserstein crossover.

As numerical examples, we considered three topology optimization problems: a two-dimensional stress minimization problem, a two-dimensional turbulent heat transfer problem, and a three-dimensional stress minimization problem.
In the first example, a comparison with the VAE-based crossover used in DDTD highlighted the superiority of the proposed Wasserstein crossover.
Across all examples, the proposed framework effectively handled non-differentiable or strongly multimodal evaluation functions and successfully derived physically reasonable optimized designs with a computationally feasible effort.
These results highlight the potential of Wasserstein crossover as a powerful crossover operator for topology optimization.
Thus, a natural future direction is to further investigate its applicability to more complex topology optimization problems, such as those involving geometric nonlinearity or buckling, unsteady thermofluid problems, as well as other physics-related problems.

One of the major challenges of Wasserstein crossover lies in its current limitation to problems with rectangular design domains.
Although the computation of optimal transport is theoretically possible~\cite{sejourne2022, bauer2024}, the acceleration achieved via convolutional operations becomes restricted, leading to impractically high computational costs.
To address this issue, normalization techniques such as design domain mapping~\cite{yamasaki2019} are expected to be beneficial, which will be explored in future work.
Additionally, another challenge in the optimization framework is the need to formulate a pseudo-optimization problem that can be easily solved using the density-based method.
Developing a systematic approach to construct the initial population in a problem-independent manner is a crucial topic, which will also be investigated in our subsequent studies.

\section*{CRediT authorship contribution statement}

\textbf{Taisei Kii:} Writing -- original draft, Visualization, Validation, Methodology, Investigation, Formal analysis.
\textbf{Kentaro Yaji:} Writing -- review \& editing, Methodology, Project administration, Supervision, Funding acquisition.
\textbf{Hiroshi Teramoto:} Writing -- review \& editing, Conceptualization, Methodology, Project administration.
\textbf{Kikuo Fujita:} Writing -- review \& editing, Supervision.

\section*{Declaration of competing interest}

The authors declare that they have no known competing financial interests or personal relationships that could have appeared to influence the work reported in this paper.

\section*{Acknowledgments}

This work was supported by JSPS KAKENHI, Japan, Grant Numbers 23H03799 and 24KJ1640.

\section*{Data availability}

Data will be available on request.

\bibliographystyle{elsarticle-num} 
\bibliography{references}

\begin{thebibliography}{10}
\expandafter\ifx\csname url\endcsname\relax
  \def\url#1{\texttt{#1}}\fi
\expandafter\ifx\csname urlprefix\endcsname\relax\def\urlprefix{URL }\fi
\expandafter\ifx\csname href\endcsname\relax
  \def\href#1#2{#2} \def\path#1{#1}\fi

\bibitem{bendsoe1988}
M.~P. Bends{\o}e, N.~Kikuchi, Generating optimal topologies in structural
  design using a homogenization method, Computer Methods in Applied Mechanics
  and Engineering 71~(2) (1988) 197--224.
\newblock \href {https://doi.org/10.1016/0045-7825(88)90086-2}
  {\path{doi:10.1016/0045-7825(88)90086-2}}.

\bibitem{bendsoe1989}
M.~P. Bends{\o}e, Optimal shape design as a material distribution problem,
  Structural Optimization 1 (1989) 193--202.
\newblock \href {https://doi.org/10.1007/BF01650949}
  {\path{doi:10.1007/BF01650949}}.

\bibitem{sethian2000}
J.~A. Sethian, A.~Wiegmann, Structural boundary design via level set and
  immersed interface methods, Journal of Computational Physics 163~(2) (2000)
  489--528.
\newblock \href {https://doi.org/10.1006/jcph.2000.6581}
  {\path{doi:10.1006/jcph.2000.6581}}.

\bibitem{le2010}
C.~Le, J.~Norato, T.~Bruns, C.~Ha, D.~Tortorelli, Stress-based topology
  optimization for continua, Structural and Multidisciplinary Optimization
  41~(4) (2010) 605--620.
\newblock \href {https://doi.org/10.1007/s00158-009-0440-y}
  {\path{doi:10.1007/s00158-009-0440-y}}.

\bibitem{alexandersen2020}
J.~Alexandersen, C.~S. Andreasen, A review of topology optimisation for
  fluid-based problems, Fluids 5~(1) (2020) 29.
\newblock \href {https://doi.org/10.3390/fluids5010029}
  {\path{doi:10.3390/fluids5010029}}.

\bibitem{thomas1996}
T.~B\"{a}ck, Evolutionary Algorithms in Theory and Practice: Evolution
  Strategies, Evolutionary Programming, Genetic Algorithms, Oxford University
  Press, 1996.
\newblock \href {https://doi.org/10.1093/oso/9780195099713.001.0001}
  {\path{doi:10.1093/oso/9780195099713.001.0001}}.

\bibitem{yu2010}
X.~Yu, M.~Gen, Introduction to Evolutionary Algorithms, Springer, 2010.
\newblock \href {https://doi.org/10.1007/978-1-84996-129-5}
  {\path{doi:10.1007/978-1-84996-129-5}}.

\bibitem{goldberg1989}
D.~E. Goldberg, Genetic Algorithms in Search, Optimization, and Machine
  Learning, Addison Wesley, 1989.

\bibitem{guirguis2020}
D.~Guirguis, N.~Aulig, R.~Picelli, B.~Zhu, Y.~Zhou, W.~Vicente, F.~Iorio,
  M.~Olhofer, W.~Matusik, C.~A. Coello~Coello, K.~Saitou, Evolutionary
  black-box topology optimization: Challenges and promises, IEEE Transactions
  on Evolutionary Computation 24~(4) (2020) 613--633.
\newblock \href {https://doi.org/10.1109/TEVC.2019.2954411}
  {\path{doi:10.1109/TEVC.2019.2954411}}.

\bibitem{chapman1994}
C.~D. Chapman, K.~Saitou, M.~J. Jakiela, Genetic algorithms as an approach to
  configuration and topology design, Journal of Mechanical Design 116~(4)
  (1994) 1005--1012.
\newblock \href {https://doi.org/10.1115/1.2919480}
  {\path{doi:10.1115/1.2919480}}.

\bibitem{jakiela2000}
M.~J. Jakiela, C.~Chapman, J.~Duda, A.~Adewuya, K.~Saitou, Continuum structural
  topology design with genetic algorithms, Computer Methods in Applied
  Mechanics and Engineering 186~(2--4) (2000) 339--356.
\newblock \href {https://doi.org/10.1016/S0045-7825(99)00390-4}
  {\path{doi:10.1016/S0045-7825(99)00390-4}}.

\bibitem{wang2005}
S.~Y. Wang, K.~Tai, Structural topology design optimization using genetic
  algorithms with a bit-array representation, Computer Methods in Applied
  Mechanics and Engineering 194~(36--38) (2005) 3749--3770.
\newblock \href {https://doi.org/10.1016/j.cma.2004.09.003}
  {\path{doi:10.1016/j.cma.2004.09.003}}.

\bibitem{madeira2010}
J.~F.~A. Madeira, H.~L. Pina, H.~C. Rodrigues, {GA} topology optimization using
  random keys for tree encoding of structures, Structural and Multidisciplinary
  Optimization 40 (2010) 227--240.
\newblock \href {https://doi.org/10.1007/s00158-008-0353-1}
  {\path{doi:10.1007/s00158-008-0353-1}}.

\bibitem{nimura2024}
N.~Nimura, A.~Oyama, Multiobjective evolutionary topology optimization
  algorithm using quadtree encoding, IEEE Access 12 (2024) 73839--73848.
\newblock \href {https://doi.org/10.1109/ACCESS.2024.3404594}
  {\path{doi:10.1109/ACCESS.2024.3404594}}.

\bibitem{wu2010}
C.-Y. Wu, K.-Y. Tseng, Topology optimization of structures using modified
  binary differential evolution, Structural and Multidisciplinary Optimization
  42 (2010) 939--953.
\newblock \href {https://doi.org/10.1007/s00158-010-0523-9}
  {\path{doi:10.1007/s00158-010-0523-9}}.

\bibitem{luh2011}
G.-C. Luh, C.-Y. Lin, Y.-S. Lin, A binary particle swarm optimization for
  continuum structural topology optimization, Applied Soft Computing 11~(2)
  (2011) 2833--2844.
\newblock \href {https://doi.org/10.1016/j.asoc.2010.11.013}
  {\path{doi:10.1016/j.asoc.2010.11.013}}.

\bibitem{sigmund2011}
O.~Sigmund, On the usefulness of non-gradient approaches in topology
  optimization, Structural and Multidisciplinary Optimization 43 (2011)
  589--596.
\newblock \href {https://doi.org/10.1007/s00158-011-0638-7}
  {\path{doi:10.1007/s00158-011-0638-7}}.

\bibitem{fujii2018}
G.~Fujii, M.~Takahashi, Y.~Akimoto, {CMA-ES}-based structural topology
  optimization using a level set boundary expression---{A}pplication to optical
  and carpet cloaks, Computer Methods in Applied Mechanics and Engineering 332
  (2018) 624--643.
\newblock \href {https://doi.org/10.1016/j.cma.2018.01.008}
  {\path{doi:10.1016/j.cma.2018.01.008}}.

\bibitem{tanaka2024}
S.~Tanaka, G.~Fujii, {CMA}-{ES}-based topology optimization accelerated by
  spectral level-set-boundary modeling, Computer Methods in Applied Mechanics
  and Engineering 432 (2024) 117331.
\newblock \href {https://doi.org/10.1016/j.cma.2024.117331}
  {\path{doi:10.1016/j.cma.2024.117331}}.

\bibitem{furuta2024}
K.~Furuta, Y.~Tsukuda, T.~Sasaki, N.~Ishida, K.~Izui, S.~Nishiwaki,
  S.~Watanabe, Differential evolution-based topology optimization using
  {Karhunen-Lo{\`e}ve} expansion, Engineering Optimization 57~(9) (2024)
  2419--2445.
\newblock \href {https://doi.org/10.1080/0305215X.2024.2400558}
  {\path{doi:10.1080/0305215X.2024.2400558}}.

\bibitem{yamasaki2021}
S.~Yamasaki, K.~Yaji, K.~Fujita, Data-driven topology design using a deep
  generative model, Structural and Multidisciplinary Optimization 64~(3) (2021)
  1401--1420.
\newblock \href {https://doi.org/10.1007/s00158-021-02926-y}
  {\path{doi:10.1007/s00158-021-02926-y}}.

\bibitem{woldseth2022}
R.~V. Woldseth, N.~Aage, J.~A. B{\ae}rentzen, O.~Sigmund, On the use of
  artificial nerural networks in topology optimization, Structural and
  Multidisciplinary Optimization 65~(10) (2022) 294.
\newblock \href {https://doi.org/10.1007/s00158-022-03347-1}
  {\path{doi:10.1007/s00158-022-03347-1}}.

\bibitem{kingma2013}
D.~P. Kingma, M.~Welling, Auto-encoding variational bayes, arXiv (2013).
\newblock \href {https://doi.org/10.48550/arXiv.1312.6114}
  {\path{doi:10.48550/arXiv.1312.6114}}.

\bibitem{goodfellow2014}
I.~J. Goodfellow, J.~Pouget-Abadie, M.~Mirza, B.~Xu, D.~Warde-Farley, S.~Ozair,
  A.~Courville, Y.~Bengio, Generative adversarial networks, arXiv (2014).
\newblock \href {https://doi.org/10.48550/arXiv.1406.2661}
  {\path{doi:10.48550/arXiv.1406.2661}}.

\bibitem{kii2024}
T.~Kii, K.~Yaji, K.~Fujita, Z.~Sha, C.~C. Seepersad, Latent crossover for
  data-driven multifidelity topology design, Journal of Mechanical Design
  146~(5) (2024) 051713.
\newblock \href {https://doi.org/10.1115/1.4064979}
  {\path{doi:10.1115/1.4064979}}.

\bibitem{kato2025}
M.~Kato, T.~Kii, K.~Yaji, K.~Fujita, Maximum stress minimization via
  data-driven multifidelity topology design, Journal of Mechanical Design
  147~(8) (2025) 081702.
\newblock \href {https://doi.org/10.1115/1.4067750}
  {\path{doi:10.1115/1.4067750}}.

\bibitem{yang2025pca}
J.~Yang, K.~Yaji, S.~Yamasaki, Data-driven topology design based on principal
  component analysis for {3D} structural design problems, Structural and
  Multidisciplinary Optimization 68~(5) (2025) 103.
\newblock \href {https://doi.org/10.1007/s00158-025-04036-5}
  {\path{doi:10.1007/s00158-025-04036-5}}.

\bibitem{yang2025ifl}
Y.~Yang, K.~Yaji, S.~Yamasaki, K.~Fujita, Image fragmented learning for
  data-driven topology design, Structural and Multidisciplinary Optimization
  68~(8) (2025) 148.
\newblock \href {https://doi.org/10.1007/s00158-025-04054-3}
  {\path{doi:10.1007/s00158-025-04054-3}}.

\bibitem{yaji2022}
K.~Yaji, S.~Yamasaki, K.~Fujita, Data-driven multifidelity topology design
  using a deep generative model: Application to forced convection heat transfer
  problems, Computer Methods in Applied Mechanics and Engineering 388 (2022)
  114284.
\newblock \href {https://doi.org/10.1016/j.cma.2021.114284}
  {\path{doi:10.1016/j.cma.2021.114284}}.

\bibitem{luo2025a}
J.-W. Luo, K.~Yaji, L.~Chen, W.-Q. Tao, Data-driven multi-fidelity topology
  design of fin structures for latent heat thermal energy storage, Applied
  Energy 377 (2025) 124596.
\newblock \href {https://doi.org/10.1016/j.apenergy.2024.124596}
  {\path{doi:10.1016/j.apenergy.2024.124596}}.

\bibitem{luo2025b}
J.-W. Luo, K.~Yaji, L.~Chen, W.-Q. Tao, Data-driven multifidelity topology
  design for enhancing turbulent natural convection cooling, International
  Journal of Heat and Mass Transfer 240 (2025) 126659.
\newblock \href {https://doi.org/10.1016/j.ijheatmasstransfer.2024.126659}
  {\path{doi:10.1016/j.ijheatmasstransfer.2024.126659}}.

\bibitem{urata2025}
K.~Urata, R.~Tsumoto, K.~Yaji, K.~Fujita, Data-driven morphological exploration
  and shape optimization for turbulent pipe systems, Journal of Mechanical
  Design 147~(9) (2025) 091704.
\newblock \href {https://doi.org/10.1115/1.4068984}
  {\path{doi:10.1115/1.4068984}}.

\bibitem{zhou2025}
D.~Zhou, K.~Nomura, S.~Yamasaki, Data-driven topology design for conductor
  layout problem of electromagnetic interference filter, IEEE Transactions on
  Electromagnetic Compatibility 67~(3) (2025) 872--883.
\newblock \href {https://doi.org/10.1109/TEMC.2025.3558260}
  {\path{doi:10.1109/TEMC.2025.3558260}}.

\bibitem{xu2025}
S.~Xu, H.~Kawabe, K.~Yaji, Evolutionary de-homogenization using a generative
  model for optimizing solid-porous infill structures considering the stress
  concentration issue, Materials \& Design 257 (2025) 114380.
\newblock \href {https://doi.org/10.1016/j.matdes.2025.114380}
  {\path{doi:10.1016/j.matdes.2025.114380}}.

\bibitem{kawabe2025}
H.~Kawabe, K.~Yaji, Y.~Aoki, Data-driven multifidelity topology design with
  multi-channel variational auto-encoder for concurrent optimization of
  multiple design variable fields, Computer Methods in Applied Mechanics and
  Engineering 437 (2025) 117772.
\newblock \href {https://doi.org/10.1016/j.cma.2025.117772}
  {\path{doi:10.1016/j.cma.2025.117772}}.

\bibitem{solomon2015}
J.~Solomon, F.~de~Goes, G.~Peyr^^c3^^a9, M.~Cuturi, A.~Butscher, A.~Nguyen,
  T.~Du, L.~Guibas, Convolutional {Wasserstein} distances: Efficient optimal
  transportation on geometric domains, ACM Transactions on Graphics 34~(4)
  (2015) 1--11.
\newblock \href {https://doi.org/10.1145/2766963} {\path{doi:10.1145/2766963}}.

\bibitem{villani2009}
C.~Villani, Optimal Transport: Old and New, Springer, 2009.
\newblock \href {https://doi.org/10.1007/978-3-540-71050-9}
  {\path{doi:10.1007/978-3-540-71050-9}}.

\bibitem{rubner2000}
Y.~Rubner, C.~Tomasi, L.~J. Guibas, The earth mover's distance as a metric for
  image retrieval, International Journal of Computer Vision 40~(2) (2000)
  99--121.
\newblock \href {https://doi.org/10.1023/A:1026543900054}
  {\path{doi:10.1023/A:1026543900054}}.

\bibitem{peyre2019}
G.~Peyr\'{e}, M.~Cuturi, Computational Optimal Transport: With Applications to
  Data Science, Now Fundations and Trends, 2019.
\newblock \href {https://doi.org/10.1561/2200000073}
  {\path{doi:10.1561/2200000073}}.

\bibitem{yaji2020}
K.~Yaji, S.~Yamasaki, K.~Fujita, Multifidelity design guided by topology
  optimization, Structural and Multidisciplinary Optimization 61~(3) (2020)
  1071--1085.
\newblock \href {https://doi.org/10.1007/s00158-019-02406-4}
  {\path{doi:10.1007/s00158-019-02406-4}}.

\bibitem{cuturi2013}
M.~Cuturi, Sinkhorn distances: Lightspeed computation of optimal transport,
  Advances in Neural Information Processing Systems 26 (2013) 2292--2300.
\newblock \href {https://doi.org/10.48550/arXiv.1406.2661}
  {\path{doi:10.48550/arXiv.1406.2661}}.

\bibitem{sinkhorn1967}
R.~Sinkhorn, Diagonal equivalence to matrices with prescribed row and column
  sums, The American Mathematical Monthly 74~(4) (1967) 402--405.
\newblock \href {https://doi.org/10.2307/2314570} {\path{doi:10.2307/2314570}}.

\bibitem{Ryu2018}
E.~K. Ryu, Y.~Chen, W.~Li, S.~Osher, Vector and matrix optimal mass transport:
  Theory, algorithm, and applications, SIAM Journal on Scientific Computing
  40~(5) (2018).
\newblock \href {https://doi.org/10.1137/17M1163396}
  {\path{doi:10.1137/17M1163396}}.

\bibitem{arjovsky2017}
M.~Arjovsky, S.~Chintala, L.~Bottou, Wasserstein {GAN}, arXiv (2017).
\newblock \href {https://doi.org/10.48550/arXiv.1701.07875}
  {\path{doi:10.48550/arXiv.1701.07875}}.

\bibitem{kusner2015}
M.~J. Kusner, Y.~Sun, N.~I. Kolkin, K.~Q. Weinberger, From word embeddings to
  document distances, in: Proceedings of the 32nd International Conference on
  Machine Learning, Vol.~37, 2015, pp. 957--966.
\newblock \href {https://doi.org/10.5555/3045118.3045221}
  {\path{doi:10.5555/3045118.3045221}}.

\bibitem{agueh2011}
M.~Agueh, G.~Carlier, Barycenters in the {Wasserstein} space, SIAM Journal on
  Mathematical Analysis 43~(2) (2011) 904--924.
\newblock \href {https://doi.org/10.1137/100805741}
  {\path{doi:10.1137/100805741}}.

\bibitem{benamou2015}
J.-D. Benamou, G.~Carlier, M.~Cuturi, L.~Nenna, G.~Peyr\'{e}, Iterative bregman
  projections for regularized transportation problems, SIAM Journal on
  Scientific Computing 37~(2) (2015) A1111--A1138.
\newblock \href {https://doi.org/10.1137/141000439}
  {\path{doi:10.1137/141000439}}.

\bibitem{rawat2019a}
S.~Rawat, M.-H.~H. Shen, A novel topology optimization approach using
  conditional deep learning, arXiv (2019).
\newblock \href {https://doi.org/10.48550/arXiv.1901.04859}
  {\path{doi:10.48550/arXiv.1901.04859}}.

\bibitem{rawat2019b}
S.~Rawat, M.-H.~H. Shen, Application of adversarial networks for {3D}
  structural topology optimization, SAE Technical Paper (2019).
\newblock \href {https://doi.org/10.4271/2019-01-0829}
  {\path{doi:10.4271/2019-01-0829}}.

\bibitem{wang2023}
Z.~Wang, S.~Melkote, D.~W. Rosen, Generative design by embedding topology
  optimization into conditional generative adversarial network, Journal of
  Mechanical Design 145~(11) (2023) 111702.
\newblock \href {https://doi.org/10.1115/1.4062980}
  {\path{doi:10.1115/1.4062980}}.

\bibitem{pereira2024}
L.~Pereira, L.~Driemeier, Wasserstein generative adversarial networks for
  topology optimization, Structures 67 (2024) 106924.
\newblock \href {https://doi.org/10.1016/j.istruc.2024.106924}
  {\path{doi:10.1016/j.istruc.2024.106924}}.

\bibitem{zeng2024}
Q.~Zeng, X.~Liu, X.~Zhu, X.~Zhang, P.~Hu, Data-driven structural topology
  optimization method using conditional {Wasserstein} generative adversarial
  networks with gradient penalty, Computer Modeling in Engineering and Sciences
  141~(3) (2024) 2065--2085.
\newblock \href {https://doi.org/10.32604/cmes.2024.052620}
  {\path{doi:10.32604/cmes.2024.052620}}.

\bibitem{sow2024}
B.~Sow, R.~Le~Riche, J.~Pelamatti, M.~Keller, S.~Zannane, Wasserstein-based
  evolutionary operators for optimizing sets of points: Application to
  wind-farm layout design, Applied Sciences 14~(17) (2024) 7916.
\newblock \href {https://doi.org/10.3390/app14177916}
  {\path{doi:10.3390/app14177916}}.

\bibitem{ma2021}
P.~Ma, T.~Du, J.~Z. Zhang, K.~Wu, A.~Spielberg, R.~K. Katzschmann, W.~Matusik,
  {DiffAqua}: A differentiable computational design pipeline for soft
  underwater swimmers with shape interpolation, ACM Transactions on Graphics
  40~(4) (2021) 1--14.
\newblock \href {https://doi.org/10.1145/3450626.3459832}
  {\path{doi:10.1145/3450626.3459832}}.

\bibitem{zadeh1963}
L.~Zadeh, Optimality and non-scalar-valued performance criteria, IEEE
  Transactions on Automatic Control 8~(1) (1963) 59--60.
\newblock \href {https://doi.org/10.1109/TAC.1963.1105511}
  {\path{doi:10.1109/TAC.1963.1105511}}.

\bibitem{haimes1971}
Y.~Haimes, On a bicriterion formulation of the problems of integrated system
  identification and system optimization, IEEE Transactions on Systems, Man,
  and Cybernetics 1~(3) (1971) 296--297.
\newblock \href {https://doi.org/10.1109/TSMC.1971.4308298}
  {\path{doi:10.1109/TSMC.1971.4308298}}.

\bibitem{kii2025}
T.~Kii, K.~Yaji, H.~Teramoto, K.~Fujita, Data-driven topology design with
  persistent homology for enhancing population diversity, International Journal
  of Mechanical Sciences 301 (2025) 110493.
\newblock \href {https://doi.org/10.1016/j.ijmecsci.2025.110493}
  {\path{doi:10.1016/j.ijmecsci.2025.110493}}.

\bibitem{deb2002}
K.~Deb, A.~Pratap, S.~Agarwal, T.~Meyarivan, A fast and elitist multiobjective
  genetic algorithm: {NSGA-II}, IEEE Transactions on Evolutionary Computation
  6~(2) (2002) 182--197.
\newblock \href {https://doi.org/10.1109/4235.996017}
  {\path{doi:10.1109/4235.996017}}.

\bibitem{verma2021}
S.~Verma, M.~Pant, V.~Snasel, A comprehensive review on {NSGA-II} for
  multi-objective combinatorial optimization problems, IEEE Access 9 (2021)
  57757--57791.
\newblock \href {https://doi.org/10.1109/ACCESS.2021.3070634}
  {\path{doi:10.1109/ACCESS.2021.3070634}}.

\bibitem{edelsbrunner2002}
H.~Edelsbrunner, D.~Letscher, A.~Zomorodian, Topological persistence and
  simplification, Discrete and Computational Geometry 28~(4) (2002) 511--533.
\newblock \href {https://doi.org/10.1007/s00454-002-2885-2}
  {\path{doi:10.1007/s00454-002-2885-2}}.

\bibitem{zomorodian2005}
A.~Zomorodian, G.~Carlsson, Computing persistent homology, Discrete and
  Computational Geometry 33~(2) (2005) 249--274.
\newblock \href {https://doi.org/10.1007/s00454-004-1146-y}
  {\path{doi:10.1007/s00454-004-1146-y}}.

\bibitem{mileyko2011}
Y.~Mileyko, S.~Mukherjee, J.~Harer, Probability measures on the space of
  persistence diagrams, Inverse Problems 27~(12) (2011) 124007.
\newblock \href {https://doi.org/10.1088/0266-5611/27/12/124007}
  {\path{doi:10.1088/0266-5611/27/12/124007}}.

\bibitem{shang2020}
K.~Shang, H.~Ishibuchi, L.~He, L.~M. Pang, A survey on the hypervolume
  indicator in evolutionary multiobjective optimization, IEEE Transactions on
  Evolutionary Computation 25~(1) (2020) 1--20.
\newblock \href {https://doi.org/10.1109/TEVC.2020.3013290}
  {\path{doi:10.1109/TEVC.2020.3013290}}.

\bibitem{bruns2001}
T.~E. Bruns, D.~A. Tortorelli, Topology optimization of non-linear elastic
  structures and compliant mechanisms, Computer Methods in Applied Mechanics
  and Engineering 190~(26--27) (2001) 3443--3459.
\newblock \href {https://doi.org/10.1016/S0045-7825(00)00278-4}
  {\path{doi:10.1016/S0045-7825(00)00278-4}}.

\bibitem{bourdin2001}
B.~Bourdin, Filters in topology optimization, International Journal for
  Numerical Methods in Engineering 50~(9) (2001) 2143--2158.
\newblock \href {https://doi.org/10.1002/nme.116} {\path{doi:10.1002/nme.116}}.

\bibitem{svanberg1987}
K.~Svanberg, The method of moving asymptotes―--a new method for structural
  optimization, International Journal for Numerical Methods in Engineering
  24~(2) (1987) 359--373.
\newblock \href {https://doi.org/10.1002/nme.1620240207}
  {\path{doi:10.1002/nme.1620240207}}.

\bibitem{lazarov2011}
B.~S. Lazarov, O.~Sigmund, Filters in topology optimization based on
  {Helmholtz}-type differential equations, International Journal for Numerical
  Methods in Engineering 86~(6) (2011) 765--781.
\newblock \href {https://doi.org/10.1002/nme.3072}
  {\path{doi:10.1002/nme.3072}}.

\bibitem{bradbury2018}
J.~Bradbury, R.~Frostig, P.~Hawkins, M.~J. Johnson, C.~Leary, D.~Maclaurin,
  G.~Necula, A.~Paszke, J.~Vander{P}las, S.~Wanderman-{M}ilne, Q.~Zhang,
  \href{http://github.com/jax-ml/jax}{{JAX}: Composable transformations of
  {P}ython+{N}um{P}y programs} (2018).
\newline\urlprefix\url{http://github.com/jax-ml/jax}

\bibitem{flamary2021}
R.~Flamary, N.~Courty, A.~Gramfort, M.~Z. Alaya, A.~Boisbunon, S.~Chambon,
  L.~Chapel, A.~Corenflos, K.~Fatras, N.~Fournier, L.~Gautheron, N.~T. Gayraud,
  H.~Janati, A.~Rakotomamonjy, I.~Redko, A.~Rolet, A.~Schutz, V.~Seguy, D.~J.
  Sutherland, R.~Tavenard, A.~Tong, T.~Vayer,
  \href{http://jmlr.org/papers/v22/20-451.html}{{POT}: Python optimal
  transport}, Journal of Machine Learning Research 22~(78) (2021) 1--8.
\newline\urlprefix\url{http://jmlr.org/papers/v22/20-451.html}

\bibitem{emmendoerfer2014}
H.~Emmendoerfer~Jr., E.~A. Fancello, A level set approach for topology
  optimization with local stress constraints, International Journal for
  Numerical Methods in Engineering 99~(2) (2014) 129--156.
\newblock \href {https://doi.org/10.1002/nme.4676}
  {\path{doi:10.1002/nme.4676}}.

\bibitem{emmendoerfer2016}
H.~Emmendoerfer~Jr., E.~A. Fancello, Topology optimization with local stress
  constraint based on level set evolution via reaction-diffusion, Computer
  Methods in Applied Mechanics and Engineering 305 (2016) 62--88.
\newblock \href {https://doi.org/10.1016/j.cma.2016.02.024}
  {\path{doi:10.1016/j.cma.2016.02.024}}.

\bibitem{chu2018}
S.~Chu, L.~Gao, M.~Xiao, Z.~Luo, H.~Li, X.~Gui, A new method based on adaptive
  volume constraint and stress penalty for stress-constrained topology
  optimization, Structural and Multidisciplinary Optimization 57~(3) (2018)
  1163--1185.
\newblock \href {https://doi.org/10.1007/s00158-017-1803-4}
  {\path{doi:10.1007/s00158-017-1803-4}}.

\bibitem{giraldo-londono2021}
O.~Giraldo-Londo\~{n}o, G.~H. Paulino, {PolyStress}: a {Matlab} implementation
  for local stress-constrained topology optimization using the augmented
  lagrangian method, Structural and Multidisciplinary Optimization 63~(4)
  (2021) 2065--2097.
\newblock \href {https://doi.org/10.1007/s00158-020-02760-8}
  {\path{doi:10.1007/s00158-020-02760-8}}.

\bibitem{norato2022}
J.~A. Norato, H.~A. Smith, J.~D. Deaton, R.~M. Kolonay, A
  maximum-rectifier-function approach to stress-constrained topology
  optimization, Structural and Multidisciplinary Optimization 65~(10) (2022)
  286.
\newblock \href {https://doi.org/10.1007/s00158-022-03357-z}
  {\path{doi:10.1007/s00158-022-03357-z}}.

\bibitem{bendsoe2003}
M.~P. Bends{\o}e, O.~Sigmund, Topology Optimization: Theory, Methods, and
  Applications, Springer, 2003.

\bibitem{kingma2014}
D.~P. Kingma, J.~Ba, Adam: A method for stochastic optimization, arXiv (2014).
\newblock \href {https://doi.org/10.48550/arXiv.1412.6980}
  {\path{doi:10.48550/arXiv.1412.6980}}.

\bibitem{wang2011}
F.~Wang, B.~S. Lazarov, O.~Sigmund, On projection methods, convergence and
  robust formulations in topology optimization, Structural and
  Multidisciplinary Optimization 43~(6) (2011) 767--784.
\newblock \href {https://doi.org/10.1007/s00158-010-0602-y}
  {\path{doi:10.1007/s00158-010-0602-y}}.

\bibitem{borrvall2003}
T.~Borrvall, J.~Petersson, Topology optimization of fluids in {Stokes} flow,
  International Journal for Numerical Methods in Fluids 41~(1) (2003) 77--107.
\newblock \href {https://doi.org/10.1002/fld.426} {\path{doi:10.1002/fld.426}}.

\bibitem{sejourne2022}
T.~S\'{e}journ\'{e}, G.~Peyr\'{e}, F.-X. Vialard, Unbalanced optimal transport,
  from theory to numerics, arXiv (2022).
\newblock \href {https://doi.org/10.48550/arXiv.2211.08775}
  {\path{doi:10.48550/arXiv.2211.08775}}.

\bibitem{bauer2024}
M.~Bauer, N.~Charon, T.~Needham, M.~Nishino, Path constrained unbalanced
  optimal transport, arXiv (2024).
\newblock \href {https://doi.org/10.48550/arXiv.2402.15860}
  {\path{doi:10.48550/arXiv.2402.15860}}.

\bibitem{yamasaki2019}
S.~Yamasaki, K.~Yaji, K.~Fujita, Knowledge discovery in databases for
  determining formulation in topology optimization, Structural and
  Multidisciplinary Optimization 59~(2) (2019) 595--611.
\newblock \href {https://doi.org/10.1007/s00158-018-2086-0}
  {\path{doi:10.1007/s00158-018-2086-0}}.

\end{thebibliography}

\end{document}